\documentclass[a4paper, 12pt]{amsart}
\usepackage[utf8]{inputenc}
\usepackage[includeheadfoot, margin=2.54cm]{geometry}
\usepackage{setspace}
\usepackage{amsmath}
\usepackage{amssymb}
\usepackage{amsthm}
\usepackage{bbm}
\usepackage{xcolor}
\usepackage{hyperref}
\usepackage{mathtools}
\usepackage{enumitem}
\usepackage{graphicx}
\usepackage{caption}
\usepackage{subcaption}
\usepackage{changepage}

\usepackage[
    giveninits=true,
    maxbibnames=99,
    sorting=anyt
]{biblatex}
\newcommand{\vecsp}[2]{\mathbb{F}_{#1}^{#2}}
\newcommand{\id}[1]{\mathbbm{1}_{#1}}
\newcommand{\st}{:\mathrm{ }}
\newcommand{\twr}{\mathrm{\mathbf{twr}}}
\newcommand{\wwz}{\mathrm{\mathbf{wwz}}}
\newcommand{\bracketed}[1]{{(#1)}}
\newcommand{\vecline}[1]{\underline{#1}}
\newcommand{\rank}{\textbf{rank}}

\newcommand{\cmnt}[1]{\textcolor{blue}{[#1]}}

\newcommand{\Img}[1]{\mathrm{Im}(#1)}

\newcommand{\codim}{\mathrm{codim}}

\newcommand{\At}{\mathrm{At}}
\newcommand{\polys}[1]{\mathcal{P}_{#1}}

\DeclareMathOperator*{\expct}{\mathbb{E}}

\DeclarePairedDelimiter\gen{\langle}{\rangle}
\DeclarePairedDelimiter\abs{|}{|}
\DeclarePairedDelimiter\norm{\|}{\|}
\DeclarePairedDelimiter\babs{\big|}{\big|}
\DeclarePairedDelimiter\set{\{}{\}}

\DeclarePairedDelimiter\ceil{\lceil}{\rceil}

\title{A note on lower bounds for arithmetic regularity partitions}
\author{V. Gladkova}

\makeatletter
    \newtheorem*{rep@theorem}{\rep@title}
    \newcommand{\newreptheorem}[2]{
    \newenvironment{rep#1}[1]{
     \def\rep@title{#2 \ref{##1}}
     \begin{rep@theorem}}
     {\end{rep@theorem}}}
\makeatother

\makeatletter
    \def\namedlabel#1#2{\begingroup
       \def\@currentlabel{#2}
       \label{#1}\endgroup
    }
\makeatother

\newtheorem{theorem}{Theorem}[section]
\newtheorem{definition}[theorem]{Definition}

\newtheorem*{mainA}{Theorem A}
\newtheorem*{mainB}{Theorem B}

\newtheorem{lemma}[theorem]{Lemma}
\newreptheorem{lemma}{Lemma}
\newtheorem*{nclaim}{Claim}
\newtheorem{claim}[theorem]{Claim}
\newtheorem{corollary}[theorem]{Corollary}
\newtheorem{prop}[theorem]{Proposition}

\theoremstyle{definition}

\begin{document}

\begin{abstract}
    This paper establishes lower bounds for two kinds of arithmetic regularity partitions, building on constructions of Green \cite{tower-type} and Hosseini, Lovett, Moshkovitz, and Shapira \cite{arl-lowerbound}. The first kind occurs in the so-called strong arithmetic regularity lemma due to Bhattcharrya, Fischer, and Lovett \cite[Theorem 4.9]{decomposition-lemmas}, which is an arithmetic analogue of the strong regularity lemma for graphs developed by  Alon, Fischer, Krivelevich, and Szegedy \cite{induced-graph}. Conlon and Fox \cite{conlon-fox}, as well as Kalyanasundaram and Shapira \cite{wowzer-srl-2}, demonstrated that there are graphs for which any strong regularity partition must have size at least a wowzer-type function in the pseudorandomness parameter, and the primary aim of this paper is to match this bound in the setting of vector spaces over finite fields. The second kind of arithmetic regularity partition originates from higher-order arithmetic regularity lemmas. The upper bounds on the size of these partitions are known to be of tower-type growth. Previous work \cite{tower-type, arl-lowerbound} demonstrated that this is unavoidable for the `linear' arithmetic regularity lemma of Green \cite{tower-type}, and the second contribution of this paper confirms that this continues to be necessary in the higher-order setting.
\end{abstract}

\maketitle

\section{Introduction}

A standard graph-theoretic argument shows that if a graph $G$ has few copies of a given subgraph $H$, then $G$ can be made completely $H$-free by only removing a small proportion of its edges \cite{removal-lemmas}. This is known as a \textit{graph removal lemma}, and the key ingredient in the standard proofs is Szemer\'edi's regularity lemma \cite{szemeredi} (but see also \cite{fox-new-removal}). 

In informal terms, for a given graph $G$, the regularity lemma provides a vertex partition such that the edges of $G$ between most pairs of vertex classes behave `pseudorandomly'. The size of such a partition only depends on the pseudorandonmness parameter $\epsilon$ but can be as large as a tower in $\epsilon^{-1}$, as shown by Gowers \cite{szemeredi-tower} (here the tower function $\twr: \mathbb{N} \rightarrow \mathbb{N}$ is defined by $\twr(0)=1$ and $\twr(i+1) = 2^{\twr(i)}$).

A related result is the \textit{induced removal lemma}, which states that if $H$ occurs with small muliplicity in $G$ as an induced subgraph, then $G$ can be made free of any induced copies of $H$ by flipping only a small proportion of the edges. This turned out to be a harder problem, and the first proof, due to Alon, Fischer, Krivelevich, and Szegedy \cite{induced-graph}, led to a strong version of the regularity lemma.

This strong regularity lemma \cite[Lemma 4.1]{induced-graph} gives two nested vertex partitions of $G$ with certain regularity properties. In particular, these properties allow us to find a subclass inside each vertex class of the coarser partition such that
\begin{itemize}
    \item the edges between \emph{all} pairs of such subclasses behave `pseudorandomly';

    \item for most pairs of subclasses, the edge density between them is close to that between the vertex classes containing them.
\end{itemize}
As a result, it becomes possible to restrict one's attention to the chosen subclasses and enjoy the pseudorandom behaviour of the edges while not losing too much information about $G$ as a whole. This allowed Alon, Fischer, Krivelevich, and Szegedy \cite{induced-graph} to prove the induced removal lemma for graphs, for which Szemer\'edi's regularity lemma appears insufficient. However, this gain comes at a considerable cost to bounds in applications: vertex partitions produced by Szemer\'edi's regularity lemma have size at worst tower-type in the inverse of the pseudorandommess parameter \cite{szemeredi-tower}, whereas Conlon and Fox \cite{conlon-fox}, as well as Kalyanasundaram and Shapira \cite{wowzer-srl-2}, showed that for the strong regularity lemma this can grow as fast as a wowzer function (in this paper, the wowzer function $\wwz:\mathbb{N} \rightarrow \mathbb{N}$ is defined recursively by $\wwz(1) = 2$ and $\wwz(i+1) = \twr(\wwz(i))$).

All these graph-theoretic results have arithmetic analogues in the setting of finite abelian groups. The first analogue of Szemer\'edi's regularity lemma was stated and proved by Green \cite{tower-type}. In a vector space over a fixed finite field $\vecsp{p}{}$, the equivalent of a graph is a subset of $\vecsp{p}{n}$ or, even more generally, a function $f: \vecsp{p}{n} \rightarrow [0,1]$; instead of a vertex partition of a graph, the arithmetic regularity lemma yields a partition of $\vecsp{p}{n}$ into cosets of a subspace such that $f$ is \textit{Fourier-uniform} on most of the cosets.

\begin{definition}[Fourier uniformity]
    \label{def:fourier}
    Let $H$ be a subspace of $\vecsp{p}{n}$. Given a function $F: \vecsp{p}{n} \rightarrow \mathbb{C}$ and elements $c, r \in \vecsp{p}{n}$, the \emph{Fourier transform of $F$ on $H+c$ at $r$} is defined as
        $$\widehat{F{|_{H+c}}}(r) = \expct_{x \in H+c} F(x) e_p(r^T x),$$
    where  $e_p(\cdot)$ denotes $\exp(2\pi i \cdot/p)$.

    A function $f: \vecsp{p}{n}\rightarrow \mathbb{C}$ is said to be \emph{$\epsilon$-uniform on $H+c$} if $\abs{\widehat{F{|_{H+c}}}(r)} \leq \epsilon$ for all $r \in \vecsp{p}{n}$, where $F = f - \expct_{x \in H+c} f(x)$.
\end{definition}

\begin{definition}[Partition regularity]
\label{def:regularity}
    Given a subspace $H$ of $\vecsp{p}{n}$, let $\mathcal{P}(H)$ denote the partition of $\vecsp{p}{n}$ into cosets of $H$. The partition $\mathcal{P}(H)$ is \emph{$\epsilon$-regular for $f$} if for all but an $\epsilon$-proportion of $c \in \vecsp{p}{n}$, $f$ is $\epsilon$-uniform on $H+c$.
\end{definition}

\noindent The precise statement of Green's arithmetic regularity lemma for $\vecsp{p}{n}$ is then as follows.

\begin{theorem}[Arithmetic regularity lemma
\cite{tower-type}]
\label{theorem:arl}
Fix $\epsilon > 0$. There exists $C=C_{arl}(\epsilon)$ with the following property. For any function $f:\vecsp{p}{n} \rightarrow [0,1]$ and subspace $H_0 \leqslant \vecsp{p}{n}$, there is a subspace $H \leqslant H_0$ of codimension at most $C$ in $H_0$ such that $\mathcal{P}(H)$ is $\epsilon$-regular for $f$.
\end{theorem}

Theorem \ref{theorem:arl} can be used to prove an \textit{arithmetic removal lemma} \cite{tower-type} where we are looking to eliminate solutions to a given system of linear equations. Specifically, the arithmetic removal lemma arising from Theorem \ref{theorem:arl} concerns linear systems of \emph{true complexity 1} \cite{true-complexity}. Such systems include, for instance, the single equation $x+y+z=0$ or the equation $x-2y+z=0$ defining a 3-term arithmetic progression, but not the system $x-2y=z$, $y-2z+w=0$ defining a 4-term arithmetic progression.

An \emph{induced arithmetic removal lemma} corresponds to removal of solutions to linear systems under specified colourings: for example, given a 3-colouring of $\vecsp{p}{n}$, we might wish to eliminate all `rainbow' 3-term arithmetic progressions, i.e.~ones in which each term has a different colour. Induced arithmetic removal lemmas for translation-invariant systems were developed in the work of Bhattacharyya, Grigorescu, and Shapira \cite{induced-over-f2}, and Bhattacharyya, Fischer, Hatami, Hatami, and Lovett \cite{trans-invariant}. This was subsequently extended to all linear systems of complexity 1 by Fox, Tidor, and Zhao \cite{induced-1}, then systems of any complexity by Tidor and Zhao \cite{full-induced}, with the caveat that certain `non-generic' solutions might be left behind; recent work of the author \cite{induced-partition-regular} shows that no such exceptions need be made when the linear system in question is partition-regular.

A key tool in all existing proofs is an arithmetic analogue of the strong regularity lemma. Where linear systems of complexity 1 are concerned, the appropriate version of strong regularity \cite[Theorem 5.4]{induced-1} implies that for any $\epsilon > 0$ and a function $f:\vecsp{p}{n} \rightarrow [0,1]$, there are two nested subspaces $W_2 \leqslant W_1$ of $\vecsp{p}{n}$ such that
\begin{itemize}
    \item inside each coset of $W_1$, there is a coset of $W_2$ that is $\epsilon$-regular for $f$;
    \item the density of $f$ on each such subcoset is close to the density on the corresponding coset of $W_1$.
\end{itemize}
In fact, the strong arithmetic regularity lemma itself asserts something more general, namely that, given a \emph{function} $\epsilon:\mathbb{N} \rightarrow (0,1)$, there are subspaces $W_2 \leqslant W_1$ of codimensions $C_2$ and $C_1$ respectively, such that $\mathcal{P}(W_2)$ is $\epsilon(C_1)$-regular for $f$, and the \emph{energies} of $\mathcal{P}(W_1)$ and $\mathcal{P}(W_2)$ are close. In the following definition, given a function $f:\vecsp{p}{n} \rightarrow [0,1]$ and a partition $\mathcal{P}$ of $\vecsp{p}{n}$, let $\expct_{}(f|\mathcal{P})(x)$ denote the density of $f$ on the unique part $P_x$ of $\mathcal{P}$ containing $x$, i.e.~$\expct_{}(f|\mathcal{P})(x) = \expct_{y \in P_x} f(y)$.

\begin{definition}[Energy]
\label{def:energy}
Let $f:\vecsp{p}{n} \rightarrow [0,1]$ be a function and $\mathcal{P}$ a partition of $\vecsp{p}{n}$. The \emph{energy of $\mathcal{P}$ with respect to $f$} is defined as $\mathcal{E}(\mathcal{P}) = \norm{\expct_{}(f|\mathcal{P})}_{L_2}^2$. For a subspace $H$ of $\vecsp{p}{n}$, the \emph{energy of $H$ with respect to $f$} is given by $\mathcal{E}(H) = \mathcal{E}(\mathcal{P}(H))$.
\end{definition}

\begin{theorem}[Strong arithmetic regularity lemma \cite{induced-1}]
\label{thm:sarl}
Fix a prime $p$, $\delta > 0$, and a non-increasing function $\epsilon: \mathbb{N} \rightarrow (0,1)$. There exist $C^\bracketed{i} = C^{\bracketed{i}}_{sarl}(p, \delta, \epsilon)$ for $i=1,2$ with the following property. For any function $f: \vecsp{p}{n} \rightarrow [0,1]$, there are subspaces $W_2 \leqslant W_1 \leqslant \vecsp{p}{n}$ of codimensions $C_2 \leq C^\bracketed{2}$ and $C_1 \leq C^\bracketed{1}$ respectively, such that
\begin{enumerate}[label=(\roman*)]
    \item $\mathcal{P}(W_2)$ is $\epsilon(C_1)$-regular for $f$;
    \item $\mathcal{E}(W_2) - \mathcal{E}(W_1) \leq \delta.$
\end{enumerate}
\end{theorem}

The work of Green \cite{tower-type} and Hosseini, Lovett, Moshkovitz and Shapira \cite{arl-lowerbound} established tower-type lower bounds on $C_{arl}(\epsilon)$ in the arithmetic regularity lemma (Theorem \ref{theorem:arl}), matching the graph-theoretic bounds \cite{szemeredi-tower} for Szemer\'edi's regularity lemma. Section \ref{section:construction} of this paper describes and builds upon the lower-bound constructions of these previous works to obtain a function that witnesses wowzer-type growth of $C^{\bracketed{1}}_{sarl}(p, \delta, \epsilon)$ in Theorem \ref{thm:sarl}, thus matching it to the growth of its graph-theoretic counterpart. Note that, since $W_2 \leqslant W_1$, the same lower bound applies to  $C^{\bracketed{2}}_{sarl}(p, \delta, \epsilon)$.

\begin{mainA}[Wowzer-type lower bound on $C_{sarl}^\bracketed{1}$]
\namedlabel{thm:sarl-intro}{A}
    Fix a prime $p$, $0 < \delta \leq 1/20p$, and a function $\epsilon: \mathbb{N} \rightarrow (0,1)$ such that $\epsilon(d) \leq \sqrt{\delta}/(80p^2(d+1))$. Then $C^{\bracketed{1}}_{sarl}(p, \delta, \epsilon) > \wwz(\lfloor \sqrt{\delta^{-1}}/10p \rfloor)$.
\end{mainA}

As noted, the arithmetic removal results described so far concern linear systems of true complexity 1, but it is possible to extend them to general linear systems. For example, an arithmetic removal lemma for systems of any complexity can by proved via hypergraph regularity lemmas \cite{rodl-hypegraph-regularity, gowers-hypergraph-regularity}, as was done by Shapira \cite{arithmetic-removal}; alternatively, one can employ so-called higher-order arithmetic regularity lemmas, developed in the work of \cite{montreal}, \cite{gowers-decompositions}, \cite{green-regularity}, \cite{gowers-wolf-higher-degree}, and \cite{gowers-wolf-quadratic}, for instance. Like Theorem \ref{theorem:arl}, these higher-order arithmetic regularity lemmas give rise to partitions of the space with certain desirable regularity properties, with such partitions referred to as \textit{arithmetic regularity partitions} in this paper. While higher-order arithmetic lemmas share many features with hypergraph regularity lemmas, there are some differences in the behaviour of regularity partitions in the arithmetic setting compared to those for hypergraphs.

For instance, consider the `quadratic' arithmetic regularity lemma \cite{montreal} and the regularity lemma for 3-uniform hypergraphs \cite{rodl-hypegraph-regularity, gowers-hypergraph-regularity}. The latter produces a partition of the vertex set, as in Szemer\'edi's regularity lemma, and additionally a partition of the set of pairs of vertices. Somewhat similarly, the quadratic arithmetic regularity partition consists of a partition of the space into cosets of a subspace, as in Theorem \ref{theorem:arl}, and additionally a quadratically-structured layer, given by simultaneous level sets of a bounded number of quadratic forms, that refines this partition further (see Definition \ref{def:quadratic-factor}). Moshkovitz and Shapira \cite{hypergraph-lower-bounds-companion, hypergraph-lower-bound} showed that the size of the vertex partition arising from 3-uniform hypergraph regularity lemmas must have wowzer-type growth (and in general, the size of the vertex partition for $k$-uniform hypergraph lemmas is a $k$-th order Ackermann function\footnote{Ackermann function of order 1 is defined as $\mathrm{Ack}_1(x) = 2^x$ and $\mathrm{Ack}_{m+1}(x)$ as $\mathrm{Ack}_m$ iterated $x$ times; in particular $\twr = \mathrm{Ack}_2$ and $\wwz = \mathrm{Ack}_3$}). By contrast, the proof of the quadratic arithmetic regularity lemma yields only a tower-type upper bound on the size of the whole partition (see Appendix \ref{appendix:upper-bounds}). Section \ref{sec:linear-layer} expands on concepts related to higher-order arithmetic regularity lemmas and contains a proof of the following result, stated here informally.

\begin{mainB}[Tower-type lower bound on the linear layer - informal]
    \namedlabel{thm:main-linear-layer}{B}
    There exists a function $f:\vecsp{p}{n} \rightarrow [0, 1]$ such that the `linear layer' of any quadratic regularity partition for $f$ must be comprised of cosets of a subspace with codimension at least tower-type in the regularity parameters.
\end{mainB}
\noindent Theorem \ref{thm:main-linear-layer} shows that the tower-type bound on the complexity of the linear layer of a quadratic regularity partition cannot in general be improved. It remains an open problem to prove a similar lower bound (or a sub-tower upper bound) on the complexity of the quadratic layer. In recent work, Terry and Wolf \cite{quadratic-cmplxt-vc2} show that under the assumption of bounded $\mathrm{VC}_2$-dimension, the size of the quadratic layer can be taken as polynomial in the regularity parameter, whereas in general it must be at least exponential. To the author's knowledge, the latter is the best lower bound currently known, although the true order of growth may well turn out to be tower-type as per the upper bound in Appendix \ref{appendix:upper-bounds}.

\subsection*{Acknowledgments} This work was supported by Harding Distinguished Postgraduate Scholars Programme and, in the later stages of preparing the manuscript, by funding associated to an Open Fellowship from the UK Engineering and Physical Sciences Research Council (EP/Z53352X/1). The author would like to thank Caroline Terry and Julia Wolf for introducing them to the problem, as well as Sean Prendiville and Julian Sahasrabudhe for several helpful comments on the first draft of this paper.

\section{A generalised lower-bound construction}
\label{section:construction}
The construction described in this section will be used to prove Theorem \ref{thm:sarl-intro}. It closely follows and builds upon the constructions of Green \cite{tower-type} and Hosseini et al.~\cite{arl-lowerbound}, but is of a more general form. With an appropriate choice of parameters, this general form allows us to recover both of the earlier constructions as well as produce a function $f$ that witnesses Theorem \ref{thm:sarl-intro}.

For $i \geq 0$, define $D_0 = 0$ and $D_{i+1} = \sum_{j=1}^{i+1} d_j$, where $d_{1} = 1$, $d_2 = 2$ and $d_{j+1} = p^{D_j-3}$ for $j \geq 2$. Observe that $D_i \geq \twr(i)$ for all $i \geq 1$. The choice of $p^{-3}$ as a multiplicative factor in this definition is motivated by the following lemma, which appears for $p=2$ in \cite[Claim 2.1]{arl-lowerbound} (or, with different constants, \cite[Lemma 10.1]{tower-type}) and is crucial to the proofs of lower bounds in \cite{tower-type} and \cite{arl-lowerbound}
\begin{lemma}
    \label{lemma:34span}
    Let $V = \vecsp{p}{d}$. There is a tuple of $p^3 d$ non-zero vectors of $V$ such that any $3/4$-proportion of them spans $V$.
\end{lemma}
\begin{proof}
The case $p=2$ is \cite[Claim 2.1]{arl-lowerbound}, so assume $p \geq 3$.
Choose non-zero vectors $v_1, \ldots, v_{p^3 d} \in V$ independently and uniformly at random. Let $U \leqslant V$ be any subspace of codimension 1, noting that each $v_i$ lies in $U$ with probability at most $1/p$. Now let $X_U$ be the random variable counting the number of $v_1, \ldots, v_{p^3d}$ that lie in $U$, so that $X_U$ is a sum of $p^3d$ Bernouilli random variables. By a Chernoff bound, e.g.~\cite[Theorem A.1.4]{probabilistic-method},
$$\mathbb{P}\left(X_U \geq 3p^3d/4 \right) \leq \exp\left(-2pd\left(3p/4-1\right)^2\right) < e^{-pd},$$
where the last inequality follows from the fact that $2(3p/4-1)^2 \geq 2\cdot 25/16 > 1$. Then by the union bound, the probability that there is such a subspace $U$ containing at least a $3/4$-proportion of the chosen vectors is at most $p^d e^{-pd} < 1$, since $e^{-p} < 1/p$ for all $p$. In other words, there is a choice of $(v_1, \ldots, v_{p^3d})$ satisfying the lemma, as required.
\end{proof}

Now let $\mathbf{e}_1, \ldots, \mathbf{e}_n$ denote the standard basis of $\vecsp{p}{n}$, and define a sequence of subspaces $\vecsp{p}{n}=H_0 \geqslant H_1 \geqslant H_2 \geqslant \ldots$ by $H_{i} = \gen{\mathbf{e}_1, \ldots, \mathbf{e}_{D_{i}}}^\perp;$ that is, each $H_i = \set{0}^{D_i} \times \vecsp{p}{n-D_i}$ and $\codim_{H_i} H_{i+1} = d_{i+1}$. Additionally, let $U_0 \leqslant U_1 \leqslant U_2 \leqslant \ldots$ be the subspaces defined by $U_i = \vecsp{p}{D_i} \times \set{0}^{n-D_i}$, so that $H_i \oplus U_i = \vecsp{p}{n}$ and each element of $U_i$ corresponds to a unique coset of $H_i$. Note that we have the following corollary of Lemma \ref{lemma:34span}.

\begin{corollary}
    \label{corollary:34span-construction}
    For each $i \geq 1$, there is a tuple $X_i = (\xi^{\bracketed{i}}_u \st u \in U_{i-1})$ of non-zero vectors $\xi^{\bracketed{i}}_u \in H_{i-1}$ such that the span of any subset $X' \subseteq X_i$ satisfying $\abs{X'} \geq 3/4 \abs{X_i}$ is equal to $\gen{\mathbf{e}_{D_{i-1}+1}, \ldots, \mathbf{e}_{D_i}}$.
\end{corollary}
\begin{proof}
    For $i \geq 2$, apply Lemma \ref{lemma:34span} to $V = \gen{\mathbf{e}_{D_{i-1}+1}, \ldots, \mathbf{e}_{D_i}} \subseteq H_{i-1}$ to obtain a set $X_i \subseteq V$ of $p^3 d_i$ vectors such that any $3/4$-proportion of them spans $V$. Since $p^3 d_i = p^{D_{i-1}} = \abs{U_{i-1}}$ by definition, these vectors can be labelled by elements of $U_{i-1}$.

    $X_1$ and $X_2$ may be defined manually. If $X_1$ is defined to be $(\mathbf{e}_1)$, then the span of $X_1$ is indeed $\gen{\mathbf{e}_{D_0+1}, \ldots, \mathbf{e}_{D_1}} = \gen{\mathbf{e}_1}$, and $X_1$ can be labelled by $U_0 = \set{\vecline{0}}$. Note that any subset of $X_1$ containing at least a $3/4$-proportion of $X_1$ must be the whole set, so $X_!$ satisfies the desired conclusion. For $i=2$, define $X_2 = (\mathbf{e}_2+m \mathbf{e}_3 \st m = 0, \ldots, p-1)$. Then the elements of $X_2$ can be labelled by $U_1 = \gen{\mathbf{e}_1}$, and any two distinct elements of $X_2$ span $\gen{\mathbf{e}_2, \mathbf{e}_3} = \gen{\mathbf{e}_{D_1+1}, \ldots, \mathbf{e}_{D_2}}$. As any $X' \subseteq X_2$ containing at least a $3/4$-proportion of $X_2$ contains $\ceil{3p/4} \geq 3$ elements, this completes the proof. $p^{d_1+d_2-3}$
\end{proof}

\begin{figure}[h]
     \centering
     \hfill
     \begin{subfigure}[t]{0.45\textwidth}
         \centering
         \includegraphics[width=0.8\textwidth]{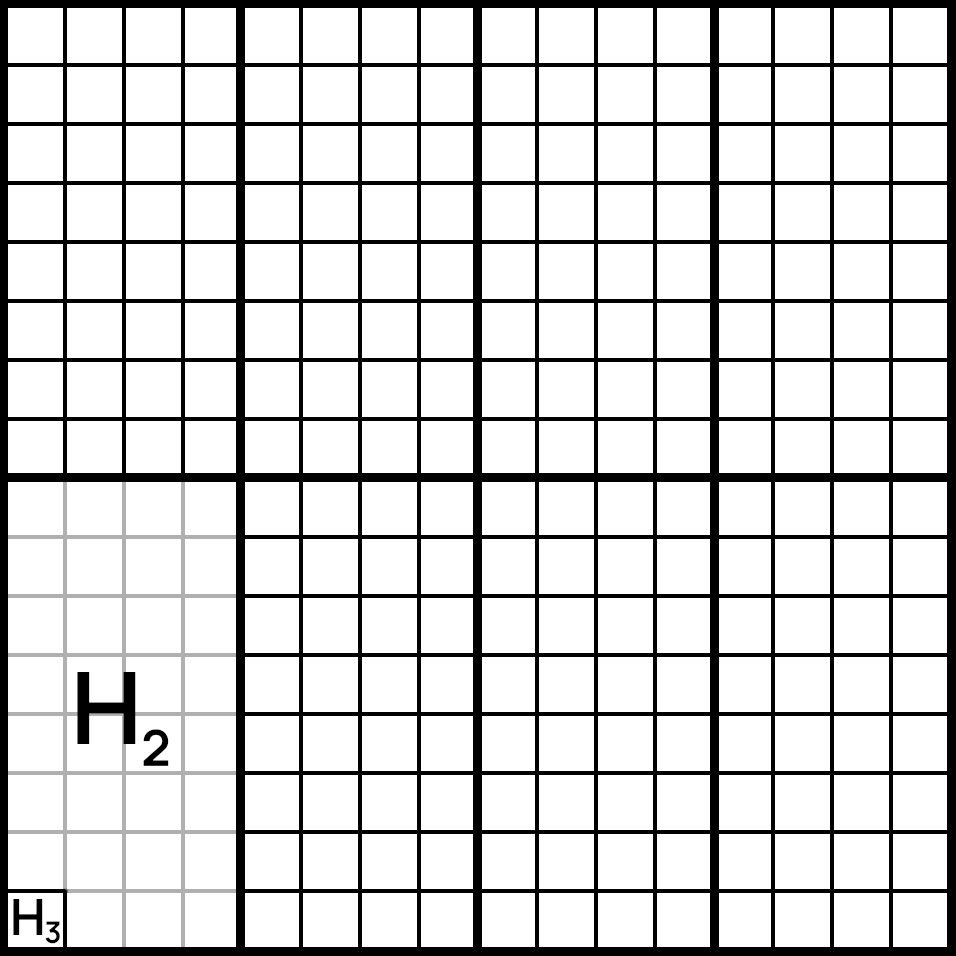}
     \end{subfigure}
     \hfill
     \begin{subfigure}[t]{0.45\textwidth}
         \centering
         \includegraphics[width=0.8\textwidth]{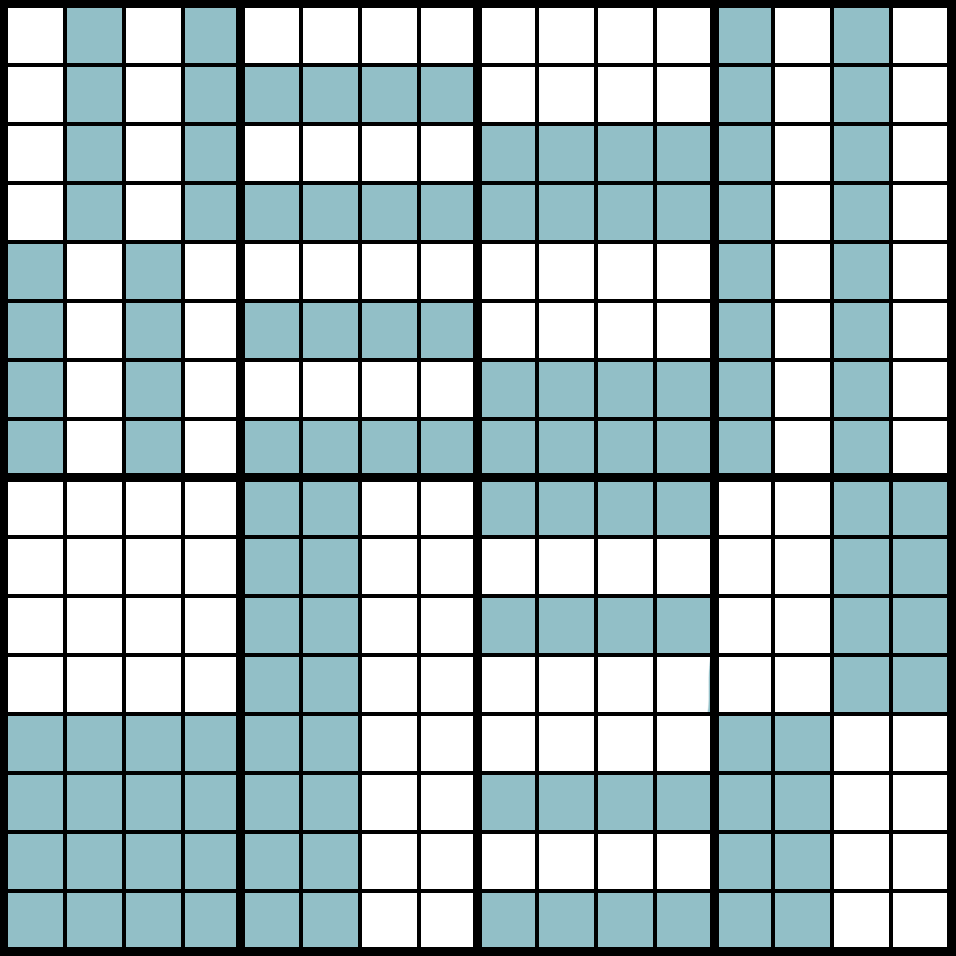}
     \end{subfigure}
     \hfill
        \caption{A possible choice of $A_3$ for the setting of $p=2$. In each coset of $H_2$, a codimension $1$ subspace is picked in such a way that $A_3$ is a union of cosets of $H_3$. (Here $H_3$ is depicted with a smaller codimension than defined, for greater visual clarity.)}
        \label{fig:sarl-layer}
\end{figure}

We will use the elements of $X_i$ to select a codimension 1 subspace in each coset of $H_{i-1}$ by taking
\begin{equation}
    \label{eq:Ai-def}
    A_i = \bigcup_{u \in U_{i-1}} \left( H_{i-1} \cap \gen{\xi^{\bracketed{i}}_u}^\perp + u \right).
\end{equation}
By the choice of $X_i$, $H_i = H_{i-1} \cap \gen{X_i}^\perp$ so, in particular, $H_i \leqslant H_{i-1} \cap \gen{\xi^{\bracketed{i}}_u}^\perp$ for all $u \in U_{i-1}$. As a result, $A_i$ is a union of cosets of $H_i$ (see Figure \ref{fig:sarl-layer}).

To complete the construction, choose some $s \in \mathbb{N}$ and weights $w_1, \ldots, w_s \in (0, 1)$ such that $\sum w_i \leq 1$, and define $f(x) = \sum_{i=1}^s w_i \id{A_i}(x)$. This function has a key property that for a given subspace $W$, $f$ is not $w_i/2p$-uniform on a large proportion of cosets of $W$ \textit{unless} $W \leqslant H_s$. This is \cite[Lemma 2.2]{arl-lowerbound}, reproved here in a more general form.

\begin{prop}
    \label{prop:large-bias}
    Let $W$ be a subspace of $\vecsp{p}{n}$ such that $W \leqslant H_{i-1}$ and $W \not\leqslant H_{i}$ for some $1 \leq i \leq s$. Then for at least a $w_{i}/8p$-proportion of $c \in \vecsp{p}{n}$, there exists a $u \in U_{i-1}$ such that $\widehat{f}|_{W+c}(\xi^\bracketed{i}_u) \geq w_{i}/2p$. In particular, $\mathcal{P}(W)$ is not $w_i/8p$-regular for $f$.
\end{prop}
\begin{proof}
    Let $S = \set{u \in U_{i-1} \st W \not\subseteq \gen{\xi^{\bracketed{i}}_u}^\perp}$ so that $S$ is the set of those $u \in U_{i-1}$ for which $W + u \not\subseteq A_i$. If $\abs{S} < \abs{X_i}/4$, then by Corollary \ref{corollary:34span-construction},
    $$\bigcap_{u \in U_{i-1} \backslash S} \gen{\xi^{\bracketed{i}}_u}^\perp = \gen{\xi^{\bracketed{i}}_u \st u \in U_{i-1} \backslash S}^\perp = \gen{\mathbf{e}_{D_{i-1}+1}, \ldots, \mathbf{e}_{D_i}}^\perp$$
    which implies $W \subseteq H_{i-1} \cap \gen{\mathbf{e}_{D_{i-1}+1}, \ldots, \mathbf{e}_{D_i}}^\perp = H_i.$ Since that would contradict the maximality of $i$, we must have $\abs{S} \geq \abs{U_i}/4$.

    Now fix some $u \in S$. Note that $W' = W \cap \gen{\xi^{\bracketed{i}}_u}^\perp$ has exactly $p$ cosets in $W$ and $\xi^{\bracketed{i}}_u \notin W^\perp$, which will be useful to keep in mind when calculating Fourier transforms. For each $h \in H_{i-1}$, let $W_{u,h}$ denote the coset $W + u + h$ of $W$ and $f_{u,h} = f|_{W_{u, h}}$ the restriction of $f$ to $W_{u,h}$.
    \begin{claim}
        $\expct_{h \in H_{i-1}} \widehat{f}_{u,h}(\xi^\bracketed{i}_u) = w_i/p$.
    \end{claim}
    \begin{proof}[Proof of Claim]
    \renewcommand{\qedsymbol}{}
        Since $f(x) = \sum_{j=1}^s w_j \id{A_j}(x)$, we can rewrite $\expct_{h \in H_{i-1}} \widehat{f}_{u,h}(\xi^\bracketed{i}_u)$ as
        \begin{equation}
        \label{eq:expct-split}
        \expct_{h \in H_{i-1}} \widehat{f}_{u,h}(\xi^\bracketed{i}_u) = \sum_{j=1}^s w_j \expct_{h \in H_{i-1}} \widehat{\id{A_j}|}_{W_{u,h}}(\xi^\bracketed{i}_u)
        \end{equation}
        and evaluate each expectation in the sum separately, depending on the value of $j$. It turns out that the only non-zero contribution comes from $j=i$.
        
        \textbf{Case 1: $j < i$.} $A_j$ is a union of cosets of $H_j$ and therefore also a union of cosets of $H_{i-1} \subseteq H_j$. As such, each coset of $H_{i-1}$ is either fully contained in $A_j$ or does not intersect it at all. Moreover, $\xi^\bracketed{i}_u \notin W^\perp$ so $\widehat{\id{A_j}|}_{W_{u,h}}(\xi^\bracketed{i}_u) = 0$ for any $h \in H_{i-1}$.

        \textbf{Case 2: $j = i$.} $A_i \cap W_{u, h} = \gen{\xi^{\bracketed{i}}_u}^\perp \cap W_{u, h}$ since $W_{u, h} \subseteq H_{i-1}+u$. On the other hand, $W' = W \cap \gen{\xi^{\bracketed{i}}_u}^\perp$ has exactly $p$ cosets in $W$, so
            $$\widehat{\id{A_j}|}_{W_{u,h}}(\xi^{\bracketed{i}}_u) = \expct_{x \in W_{u, h}} \id{\gen{\xi^{\bracketed{i}}_u}^\perp}(x) e_p\left(x^T \xi^{\bracketed{i}}_u\right) = \expct_{x \in W_{u, h}} \id{\gen{\xi^{\bracketed{i}}_u}^\perp}(x) = \frac{\abs{W'_{u, h}}}{\abs{W_{u, h}}} = \frac1p.$$

        \textbf{Case 3: $j > i$.} This is where the expectation over $h \in H_{i-1}$ plays a role as $\widehat{\id{A_j}|}_{W_{u,h}}(\xi^\bracketed{i}_u)$ may not be the same for different $h$. Firstly, observe that for any function $g$, $\expct_{h \in H_{i-1}} \expct_{x \in W_{u, h}} g(x) = \expct_{w \in W} \expct_{h \in H_{i-1}} g(w+h+u) = \expct_{x \in H_{i-1}} g(x+u).$ Then
            \begin{equation}
            \label{eq:expct-swap}
               \expct_{h \in H_{i-1}} \widehat{\id{A_j}|}_{W_{u,h}}(\xi^\bracketed{i}_u) = \expct_{x \in H_{i-1}}\id{A_j}(x+u) e_p\left({(x+u)^T\xi^{\bracketed{i}}_u}\right).
            \end{equation}
            Secondly, recall that, by definition, $A_j$ is a union over $u' \in U_{j-1}$ of the cosets given by
            $$H_{j-1} \cap \gen{\xi^\bracketed{j}_{u'}}^\perp+u' = (H_{j-1}+u') \cap \gen{\xi^\bracketed{j}_{u'}}^\perp,$$
            where the right-hand side follows from the fact that $\xi^\bracketed{j}_{u'} \in H_j$ by Corollary \ref{corollary:34span-construction} and $U_{j-1} \subseteq U_{j}$, so $u' \in \gen{\xi^\bracketed{j}_{u'}}^\perp$.
            As each $x+u$ is contained in a unique coset $H_{j-1} + u_x$ for some $u_x \in U_{j-1}$, we can write $\id{A_j}(x+u) = \id{\gen{\xi^\bracketed{j}_{u_x}}^\perp}(x+u)$. On the other hand, $U_{i-1} \subseteq U_{j-1} \subseteq \gen{\xi^\bracketed{j}_{u_x}}^\perp$, so $x+u \in \gen{\xi^\bracketed{j}_{u_x}}^\perp$ if and only if $x \in \gen{\xi^\bracketed{j}_{u_x}}^\perp$. In particular,
            $\id{A_j}(x+u) =  \expct_{\lambda \in \vecsp{p}{}} e_p\big{(}\lambda x^T \xi^\bracketed{j}_{u_x}\big{)}$. Then \eqref{eq:expct-swap} may be rewritten as
            \begin{equation}
            \label{eq:midclaim}
                \expct_{h \in H_{i-1}} \widehat{\id{A_j}|}_{W_{u,h}}(\xi^\bracketed{i}_u) = \expct_{\lambda \in \vecsp{p}{}} \expct_{x \in H_{i-1}} e_p\left(x^T (\lambda \xi^\bracketed{j}_{u_x}+\xi^\bracketed{i}_{u})\right).
            \end{equation}
            
            To evaluate this expectation, split it into a sum over cosets of $H_{j-1}$. Specifically, let $U' = U_{j-1} \cap H_{i-1}$ so that $H_{j-1} \oplus U' = H_{i-1}$, i.e.~the elements of $U'$ uniquely correspond to the cosets of $H_{j-1}$ in $H_{i-1}$. Note that for each $u' \in U'$ and $y \in H_{j-1}$, the value of $e_p\big{(}(y+u')^T\xi^\bracketed{i}_{u}\big{)}$ is equal to $c_{u'} = e_p\big{(}u'\kern0.3ex^T\xi^\bracketed{i}_{u}\big{)}$, which is independent of $y$. Therefore, since each $x \in H_{i-1}$ can be written as $x = y+u_x$ for some $y \in H_{j-1}$ and $u_x \in U'$, $e_p\big{(}x^T\xi^\bracketed{i}_{u}\big{)} = c_{u_x}$. This allows us to rewrite the right-hand side of \eqref{eq:midclaim} as
            \begin{equation*}
                \expct_{\lambda \in \vecsp{p}{}} \expct_{x \in H_{i-1}} \left[\sum_{u' \in U'}c_{u_x}\id{H_{j-1}+u'}(x)e_p\left(\lambda x^T \xi^\bracketed{j}_{u_x}\right)\right] = \expct_{u' \in U'} c_{u'} \left[ \expct_{\lambda \in \vecsp{p}{}} \expct_{x \in H_{j-1}+u'} e_p\left(\lambda x^T \xi^\bracketed{j}_{u'}\right) \right].
            \end{equation*}

            It follows from the conclusion of Corollary \ref{corollary:34span-construction} that $\xi^\bracketed{j}_{u'} \in \gen{\mathbf{e}_{D_{j-1}+1}, \ldots, \mathbf{e}_{D_j}}$ so in particular, $\xi^\bracketed{j}_{u'} \notin H_{j-1}^\perp$. As a result, $\expct_{x \in H_{j-1}+u'} e_p\big{(}\lambda x^T \xi^\bracketed{j}_{u'}\big{)}$ is non-zero if and only if $\lambda=0$, so the expectation in the square brackets on the right-hand side is equal to $1/p$. Likewise, $\xi^\bracketed{i}_{u}$ is contained in $\gen{\mathbf{e}_{D_{i-1}+1}, \ldots, \mathbf{e}_{D_i}} \subseteq \gen{\mathbf{e}_{D_{i-1}+1}, \ldots, \mathbf{e}_{D_{j-1}}} = U'$ so $\xi^\bracketed{i}_{u} \notin U'^\perp$. Hence, $\expct_{u' \in U'} c_{u'} = \expct_{u' \in U'} e_p\big{(}u'\kern0.3ex^T\xi^\bracketed{i}_{u}\big{)} = 0$, and equation \eqref{eq:midclaim} results in $\expct_{h \in H_{i-1}} \widehat{\id{A_j}|}_{W_{u,h}}(\xi^\bracketed{i}_u) = 0$.
            
            Substituting the results of each case into \eqref{eq:expct-split} proves the claim.
    \end{proof}

    \noindent To complete the proof of the proposition, let $\alpha$ denote the proportion of $h \in H_{i-1}$ for which $\widehat{f}_{u,h}(\xi^\bracketed{i}_u) \geq w_i/2p$. Then $\alpha > w_i/2p$ since
    $$\frac{w_i}p = \expct_{h \in H_{i-1}}\widehat{f}_{u,h}(\xi^\bracketed{i+1}_u) < \frac{w_i}{2p}(1-\alpha) + \alpha < \frac{w_i}{2p} + \alpha.$$
    
    This means that for each $u \in S$ there is at least a $w_i/2p$-proportion of $h \in H_{i-1}$ such that $\widehat{f}_{u,h}(\xi^\bracketed{i}_u) \geq w_i/2p$. While some $h \in H_{i-1}$ might correspond to the same coset of $W$, each such coset is counted the same number of times, so in fact for each $u \in S$, there is at least a $w_i/2p$-proportion of cosets of $W$ in $H_{i-1}+u$ on which $f$ is not Fourier-uniform. Combining this with the fact that $\abs{S} \geq \abs{U_{i-1}}/4$ gives the required $w_i/8p$-proportion of such cosets in the whole of $\vecsp{p}{n}$.
\end{proof}

\noindent When $p=2$, taking all $w_i = 16\epsilon$ and $s = \lfloor\epsilon^{-1}/16\rfloor$ recovers the construction of \cite{arl-lowerbound}, which requires tower-type growth for the codimension of $H$ in Theorem \ref{theorem:arl} since any such $H$ must satisfy $\codim\kern0.3ex H \geq \codim\kern0.3ex H_s = D_s > \twr(s-1)$ by Proposition \ref{prop:large-bias}. To obtain the same bound for $p \geq 3$, take $w_i = 8p\epsilon$ and $s = \lfloor\epsilon^{-1}/8p\rfloor$, resulting in the following.

\begin{corollary}
    \label{corr:hosseini-for-p}
    For all $\epsilon > 0$, $n \in \mathbb{N}$ and all prime $p$, there is a function $f: \vecsp{p}{n} \rightarrow [0,1]$ for which the following holds. If $H \leqslant \vecsp{p}{n}$ is a subspace such that $\mathcal{P}(H)$ is $\epsilon$-regular for $f$, then $\codim\kern0.3ex H > \twr(\lfloor\epsilon^{-1}/8p\rfloor - 1)$.
\end{corollary}

We will show that a different choice of $s$ and weights $w_i$ establishes a wowzer-type lower bound for the strong arithmetic regularity lemma. The precise choice of parameters for this purpose is inspired by the construction of Conlon and Fox \cite{conlon-fox} in the graph-theoretic setting.

Fix $0 < \delta \leq 1/20p$ and a non-increasing function $\epsilon: \mathbb{N} \rightarrow (0,1)$ such that $\epsilon(0) \leq \sqrt{\delta}/80p^2$. Let $t = \lfloor \sqrt{\delta^{-1}}/10p \rfloor$ and $I = \set{(i,j) \st 1 \leq i \leq t, 1 \leq j \leq h_i}$ for $h_1, \ldots, h_{t}$ to be picked later. The set $I$ provides a more convenient way to refer to the various parts of the construction via the translation $\phi:I \rightarrow [1,t]$ given by $\phi(i,j) = \sum_{k=1}^{i-1} {h_k} + j$. Specifically, we will write $w_{i, j} = w_{\phi(i,j)}$ and $H_{i, j} = H_{\phi(i, j)}$ so that, in particular, $H_{t,h_t} = H_s$ where $s = \sum_{k=1}^{t} {h_k}$, and the chain of subspaces looks as follows:
\begin{equation*}
\vecsp{p}{n} = H_{0} \geqslant \underbrace{H_{1,1} \geqslant \ldots \geqslant H_{1, h_1}}_{i=1} \geqslant \underbrace{H_{2,1} \geqslant \ldots \geqslant H_{2, h_2}}_{i=2} \geqslant \ldots \geqslant \underbrace{H_{t, 1} \geqslant \ldots \geqslant H_{t, h_t}}_{i=t} = H_s.
\end{equation*}
With this notation, $f$ can be rewritten as $f = \sum_{i=1}^t \sum_{j=1}^{h_i} w_{i,j} \id{A_{\phi(i,j)}}$.

Let $\epsilon_1, \ldots, \epsilon_t$ and $h_1, \ldots, h_t$ be defined by setting $\epsilon_1 = \epsilon(0)$, $h_i = \lfloor\sqrt{\delta}\epsilon_i^{-1}/8p\rfloor$ and $\epsilon_{i+1} = \epsilon(D_{\phi(i, h_i-9p)})$. Note that all $h_i$ are at least $h_1 \geq 10p$ so $\phi(i, h_i-9p)$ is well-defined.  Finally, define the weights by
    $$w_{i,j} = \begin{cases}
        8p \epsilon_i &\text{ if }j < h_i-9p,\\
        \max(8p \epsilon_i, \sqrt{\delta}) &\text{ otherwise}.
    \end{cases}$$
We can verify that these weights add up to at most $1$ as
$$\sum_{i=1}^t \sum_{j=1}^{h_i} w_{i,j} \leq \sum_{i=1}^t \Big[8p \epsilon_i h_i + 9p\sqrt{\delta}\Big] \leq  \lfloor \sqrt{\delta^{-1}}/10p \rfloor (\sqrt{\delta} + 9p\sqrt{\delta}) \leq 1,$$
so $f$ defined by such parameters is a valid instance of the construction presented in this section.

\section{Proof of Theorem \ref{thm:sarl-intro}}
\label{sec:sarl-wowzer-bound}

While Proposition \ref{prop:large-bias} was sufficient to prove a lower bound on the arithmetic regularity lemma, Theorem \ref{thm:sarl} requires some additional information regarding the energies of $\mathcal{P}(H_1), \ldots, \mathcal{P}(H_s)$. To this end, we will prove two auxiliary results before showing that the function $f$ as defined at the end of Section \ref{section:construction} requires wowzer-type codimensions in Theorem \ref{thm:sarl}. The following are standard properties of energy (see, for instance, \cite[Proposition 5.2]{induced-1}).
\begin{lemma}[Properties of energy]
    \label{lemma:energy-props}
    Let $g:\vecsp{p}{n} \rightarrow [-1,1]$ be a function, and let $\mathcal{P}$ and $\mathcal{Q}$ be partitions of $\vecsp{p}{n}$ such that $\mathcal{Q}$ refines $\mathcal{P}$. Then the energies $\mathcal{E}(\mathcal{P})$ and $\mathcal{E}(\mathcal{Q})$ satisfy
    \begin{enumerate}[label=(\roman*)]
        \item \label{energyprop:bounded} $0 \leq \mathcal{E}(\mathcal{P}) \leq 1$;
        \item \label{energyprop:energy-increase} $\mathcal{E}(\mathcal{Q}) - \mathcal{E}(\mathcal{P}) \geq 0$;
        \item \label{energyprop:pythagoras}  (Pythagoras' Theorem) $\mathcal{E}(\mathcal{Q}) - \mathcal{E}(\mathcal{P}) = \norm{\expct_{}(g|\mathcal{Q})-\expct_{}(g|\mathcal{P})}^2_{L_2}$.
    \end{enumerate}
\end{lemma}

\noindent The first auxiliary result establishes the energy gap between $\mathcal{P}(H_{i+1})$ and $\mathcal{P}(H_{i})$.
\begin{prop}
    \label{prop:energy-middle}
    For any $1 \leq i \leq s$, $\mathcal{E}(H_i) - \mathcal{E}(H_{i-1}) \geq {w_i}^2/p^2$.
\end{prop}
\begin{proof}
    Fix $u \in \vecsp{p}{n}$ and $1 \leq j \leq s$. Let $\alpha_{H_i+u}$ denote the density of $f$ on $H_i+u$ and write 
    $$\alpha_{H_i+u}(j) = w_j \frac{\abs{A_j \cap (H_i+u)}}{\abs{H_i}}$$
    so that $\alpha_{H_i+u} = \sum_{k=1}^s \alpha_{H_i+u}(k)$. By definition, each $A_j$ is a union of cosets of $H_j$ such that for any $v \in \vecsp{p}{n}$, $A_j \cap \left(H_{j-1}+v\right)$ consists of exactly $1/p$-proportion of all cosets of $H_j$ in $H_{j-1}+v$.

    If $j > i$, then $H_{j-1} \leqslant H_{i} \leqslant H_{i-1}$ so $A_j \cap \left(H_i + u\right)$ and $A_j \cap \left(H_{i-1} + u\right)$ both consist of exactly $1/p$-proportion of all cosets of $H_j$ in $H_i+u$ and $H_{i-1}+u$ respectively. Therefore  $\alpha_{H_{i-1}+u}(j) = \alpha_{H_i+u}(j) = w_j/p$.
    
    If $j < i$, then $H_{i} \leqslant H_{i-1} \leqslant H_j$ so $\alpha_{H_{i-1}+u}(j) = \alpha_{H_i+u}(j) \in \set{0, w_j}$ depending on whether $H_j+u \subseteq A_j$. By similar reasoning, $\alpha_{H_{i-1}+u}(i) = w_i/p$ and $\alpha_{H_i+u}(i) \in \set{0, w_i}$, which results in the lower bound
    $\abs{\alpha_{H_i+u} - \alpha_{H_{i-1}+u}} = \abs{\alpha_{H_i+u}(i) - \alpha_{H_{i-1}+u}(i)} \geq {w_i}/p.$
    
    Finally, by Lemma \ref{lemma:energy-props}\ref{energyprop:pythagoras}, $\mathcal{E}(H_{i+1}) - \mathcal{E}(H_i) = \expct_{u \in \vecsp{p}{n}} (\alpha_{H_i+u} - \alpha_{H_{i-1}+u})^2 \geq {w_i^2}/p^2$.
\end{proof}

\noindent The second result shows that the energy of $\mathcal{P}(H_i)$ is close to being maximal among all partitions into cosets with the same (or smaller) codimension.

\begin{prop}
    \label{prop:energy-start}
    Let $1 \leq i < s$ and suppose that for all $j > i$, $w_j \leq w_i$. If $W$ is a subspace of $\vecsp{p}{n}$ such that $\codim\kern0.3ex W \leq \codim\kern0.3ex H_i$, then $\mathcal{E}(W) < \mathcal{E}(H_i) + 8{w_i^2}$.
\end{prop}
\begin{proof}
    For $0 \leq j \leq s$, write $W_j = W \cap H_j$ and let $U_{W_j} \leqslant H_j$ be any subspace such that $W_j \oplus U_{W_j} = H_j$, so that the elements of $U_{W_j}$ correspond to the cosets of $W_j$ in $H_j$. By assumption, $\codim_{H_i} W_i \leq \codim\kern0.3ex H_i = D_i$, so there are at most $D_i$ linearly independent $\xi \in H_i$ such that $W_i \subseteq \gen{\xi}^\perp$. With the same notation as in the proof of Proposition \ref{prop:energy-middle}, Lemma \ref{lemma:energy-props}\ref{energyprop:pythagoras} gives
    \begin{equation}
        \label{eq:energy-start}
        \mathcal{E}(W_i) - \mathcal{E}(H_i) = \expct_{x \in \vecsp{p}{n}}(\alpha_{W_i+x} - \alpha_{H_i+x})^2 = \expct_{u \in U_i} \expct_{u' \in U_{W_i}} (\alpha_{W_i+u+u'} - \alpha_{H_i+u})^2.
    \end{equation}
    For any $1 \leq k, j \leq s$ and $u \in U_k$, it follows from the definition of $A_j$ in \eqref{eq:Ai-def} that
    \begin{equation*}
        \alpha_{H_k+u}(j) =
        \begin{cases}
            0 &\text{ if }j \leq k \text{ and } u \notin \gen{\xi^\bracketed{j}_{u}}^\perp;\\
            w_j &\text{ if }j \leq k \text{ and } u \in \gen{\xi^\bracketed{j}_{u}}^\perp;\\
            w_j/p &\text{ if } j > k,
        \end{cases}
    \end{equation*} 
    which implies that $\alpha_{W_i+u+u'}(j) = \alpha_{H_i+u}(j)$ whenever $j \leq i$, since $W_i+u+u' \subseteq H_i+u$ for all $u \in U_i$ and $u' \in U_{W_i}$. As a result,
    \vspace{-0.2cm}
    \begin{equation}
    \label{eq:energy-start-density-diff}
        \alpha_{W_i+u+u'} - \alpha_{H_i+u} = \sum_{j=i+1}^s (\alpha_{W_i+u+u'}(j) - w_j/p).
    \end{equation}

    \begin{claim}
        \label{claim:Wj-cosets}
        For each $j > i$, $\alpha_{W_i+u+u'}(j) - w_j/p = 0$ unless there is some $z \in U_{j-1}$ such that $H_{j-1}+z \subseteq H_i+u$ and $W_{j-1} \subseteq \gen{\xi^\bracketed{j}_{z}}^\perp$.
    \end{claim}
    \begin{proof}[Proof of Claim]
    \renewcommand{\qedsymbol}{}
        Using $v = u+u'$ as a shorthand, we can rewrite
        \begin{equation}
            \label{eq:claim-Wj-cosets}
            \alpha_{W_i+v}(j) = \sum_{v' \in V_{j-1}}\frac{\abs{W_{j-1}}}{\abs{W_i}}\alpha_{W_{j-1}+v+v'}(j),
        \end{equation}
        where $V_{j-1} \leqslant W_i$ is any subspace satisfying $W_{j-1} \oplus V_{j-1} = W_i$. Note that there is a unique $z_{v'} \in U_{j-1}$ such that $H_{j-1}+v+v' = H_{j-1}+z_{v'}$, and \eqref{eq:Ai-def} implies that
        $$A_j \cap \left(W_{j-1}+v+v'\right) = A_j \cap \left(H_{j-1}+z_{v'}\right) \cap \left(W_{j-1}+v+v'\right) = \gen{\xi^\bracketed{j}_{z_{v'}}}^\perp \cap \left(W_{j-1}+v+v'\right).$$
        In particular, if $W_{j-1} \not\subseteq \gen{\xi^\bracketed{j}_{z_{v'}}}^\perp$, then $\alpha_{W_{j-1}+v+v'}(j) = w_j/p$. As a result, \eqref{eq:claim-Wj-cosets} gives $\alpha_{W_i+u+u'}(j) - w_j/p = 0$ unless $W_{j-1} \subseteq \gen{\xi^\bracketed{j}_{z_{v'}}}^\perp$ for some $v' \in V_{j-1}$. Finally, we have $V_{j-1} \leqslant W_i \leqslant H_i$ so $H_{j-1}+z_{v'} = H_{j-1}+v+v' \subseteq H_i+u$, which proves the claim.
    \end{proof}
    
    \noindent In the cases where $\alpha_{W_i+u+u'}(j) - w_j/p \neq 0$, the bound $\abs{\alpha_{W_i+v+u}(j) - w_j/p} < w_i$ can be used instead, since $\alpha_{W_i+v+u}(j) \in [0, w_j]$ and $w_i \geq w_j$ for all $j > i$. Thus, by the triangle inequality applied to \eqref{eq:energy-start-density-diff}, $\abs{\alpha_{W_i+u+u'} - \alpha_{H_i+u}} < w_i \abs{J_{u}}$ where
    $$J_u = \set{j \st \exists z \in U_{j-1} \text{ s.t.~} H_{j-1}+z \subseteq H_i+u \text{ and } W_{j-1} \subseteq \gen{\xi^\bracketed{j}_{z}}^\perp} \cap [i+1, s].$$
    Let $F = \set{u \in U_i \st J_u \neq \emptyset}$, $J = \bigcup_{u \in U_i} J_u$ and $N = \abs{J}$. Substituting into \eqref{eq:energy-start} results in
    \begin{equation}
        \label{eq:energy-diff-prop}
        \mathcal{E}(W_i) - \mathcal{E}(H_i) < {w_i^2} \expct_{u \in U_i} \abs{J_u}^2 \leq {w_i^2} \frac{\abs{F}}{\abs{U_i}} N^2.
    \end{equation}
    In order to bound the size of $F$, we will consider a different set, one that can be tied to the codimension of $W_i$. For each $u \in F$ and $j \in J_u$, let $\xi^\bracketed{j}_{z_u}\in X_{j}$ be such that $H_{j-1}+z_u \subseteq H_i+u \text{ and } W_{j-1} \subseteq \gen{\xi^\bracketed{j}_{z_u}}^\perp$, which exists by the definition of $J_{u}$. Then define $S = \set{\xi_{z_u}^{\bracketed{j}} \st u \in F, j \in J_u}$, $F(j) = \set{u \in U_i \st j \in J_u}$, and $S(j) = \set{\xi_{z_u}^{\bracketed{j}} \st u \in F(j)}$. In particular, we have $S = \bigcup_{j \in J} S(j)$ and $F = \bigcup_{j \in J} F(j)$.
    
    If $u, u' \in F(j)$ are distinct, then $\xi_{z_u}^{\bracketed{j}} \neq \xi_{z_{u'}}^{\bracketed{j}}$ since otherwise we would have $z_u = z_{u'}$, which would imply $H_{j-1}+z_{u} \subseteq (H_i+u) \cap (H_i+u') = \emptyset$. As such, every $u \in F(j)$ corresponds to a distinct $\xi_{z_u}^{\bracketed{j}}$, so $\abs{S(j)} = \abs{F(j)}$. Additionally, $\gen{S(j)} \cap \gen{S(j')} = \set{0}$ for distinct $j, j' \in J$ since $X_{j'} \subseteq \gen{X_{j}}^\perp$ whenever $j < j'$.

    As a consequence, $\abs{F} \leq \sum_{j \in J} \abs{F(j)} = \sum_{j \in J} \abs{S(j)} = \abs{S}$. Moreover, if $m$ and $m_j$ are the maximum numbers of linearly independent vectors in $S$ and $S(j)$ respectively, then $m = \sum_{j \in J} m_j$ and $\abs{S} \leq \sum_{j \in J} p^{m_j}$.
    
    \begin{claim}
        $m \leq \codim_{H_i} W_i.$
    \end{claim}
    \begin{proof}[Proof of Claim]
    \renewcommand{\qedsymbol}{}
        We will show by induction that $\abs{W_k} \leq p^{n-D_k-\sum_{j=k+1}^s m_j}$ for any $k$ such that $i \leq k \leq s-1$. The base case $k=s-1$ holds since $W_{s-1} \leqslant H_{s-1} \cap \gen{\xi \in S(s)}^\perp$ and $\abs{H_{s-1} \cap \gen{\xi \in S(s)}^\perp} = p^{n-D_{s-1}-m_s}$.
    
        Now suppose that $\abs{W_k} \leq p^{n-D_k-\sum_{j=k+1}^s m_j}$ for some $i < k < s$. On the one hand, $H_{k} = H_{k-1} \cap \gen{\xi \in X_{k}}^\perp \leqslant H_{k-1} \cap \gen{\xi \in S(k)}^\perp$ so that in particular,
        $$H_{k-1} \cap \gen{\xi \in S(k)}^\perp = \bigcup_{c \in C_{k}} H_{k} + c$$
        for some set $C_{k} \subseteq U_{k}$. Note that $\abs{C_{k}} = \abs{H_{k-1} \cap \gen{\xi \in S(k)}^\perp}/\abs{H_{k}} = p^{D_{k}-D_{k-1}-m_{k}}$.        On the other hand, $W_{k-1} \leqslant H_{k-1} \cap \gen{\xi \in S(k)}^\perp$ so
        \begin{equation*}
            W_{k-1} = W_{k-1} \cap H_{k-1} \cap \gen{\xi \in S(k)}^\perp = \bigcup_{c \in C_{k}} W_{k} + c,
        \end{equation*}
        and therefore $\abs{W_{k-1}} = \abs{W_{k}}\abs{C_{k}} \leq p^{n-D_{k-1}-\sum_{j=k}^s m_j}$, which completes the inductive step. As a result, $\abs{W_i} \leq p^{n-D_i-\sum_{j=i+1}^s m_j} = p^{n-D_i-m}$, and $m \leq \codim_{H_i} W_i$ as required.
    \end{proof}
    
    \noindent By the initial assumption on the codimension of $W_i$, it follows that $m \leq D_i$. On the other hand, for any $k, j \in J$, $m_k \geq 1$ and $m_j = m - \sum_{k \neq j} m_k \leq m-N+1 \leq D_i-N+1$. Substituting $\abs{F} \leq \abs{S} \leq N p^{D_i-N+1}$ into \eqref{eq:energy-diff-prop} results in the energy gap of less than $w_i^2 N^3 p^{-N+1}$. It may be verified with standard calculus techniques that the function $F_p:\mathbb{R}^{+} \rightarrow \mathbb{R}^{+}$ given by $F_p(x) = x^3 p^{-x+1}$ is maximised when $p=2$ and $x=3/\log(2) \approx 4$. Therefore, $\mathcal{E}(W) - \mathcal{E}(H_i) \leq \mathcal{E}(W_i) - \mathcal{E}(H_i) < F_2(4) w_i^2 = 8 w_i^2$, as required.
\end{proof}

Let $f: \vecsp{p}{n} \rightarrow [0,1]$ be defined as at the end of Section \ref{section:construction}, and let $W_1 \leqslant W_2$ be subspaces satisfying Theorem \ref{thm:sarl} for $f$. Used in conjunction, Proposition \ref{prop:energy-middle} and Proposition \ref{prop:energy-start} allow us to deduce that if $\codim\kern0.3ex W_1 \leq \codim\kern0.3ex H_{i, h_i-9p}$ and $H_{i, h_i} \leqslant W_2$, then the energy gap $\mathcal{E}(W_2) - \mathcal{E}(W_1)$ is relatively large -- specifically, greater than $\delta$. Then $H_{i, h_i} \leqslant W_2$ can only hold if $\codim\kern0.3ex W_1 \geq \codim\kern0.3ex H_{i, h_i-9p} = D_{\phi(i, h_i-9p)}$. On the other hand, $\epsilon_{i+1} = \epsilon(D_{\phi(i, h_i-9p)})$ so if $\mathcal{P}(W_2)$ is an $\epsilon(\codim\kern0.3ex W_1)$-regular partition, it must also be $\epsilon_{i+1}$-regular, and therefore $H_{i+1, h_{i+1}} \leqslant W_2$ by Proposition \ref{prop:large-bias}. Continuing in this way leads to the bound $\codim\kern0.3ex W_1 \geq \codim\kern0.3ex H_{t-1, h_{t-1}-9p} = D_{\phi(t-1, h_{t-1}-9p)}$, as demonstrated in more detail in the proof below.

\begin{theorem}[Lower bound for strong regularity]
\label{thm:sarl-lowerbound}
Fix $0 < \delta \leq 1/20p$ and a non-increasing function $\epsilon: \mathbb{N} \rightarrow (0,1)$ such that $\epsilon(0) \leq \sqrt{\delta}/80p^2$. There exists a function $f: \vecsp{p}{n} \rightarrow [0,1]$ for which the following holds. If $W_2 \leqslant W_1 \leqslant \vecsp{p}{n}$ are subspaces satisfying
\begin{enumerate}[label=(\roman*)]
    \item $\mathcal{P}(W_2)$ is an $\epsilon(C)$-regular partition for $f$, where $C = \codim\kern0.3ex W_1$;
    \item $\mathcal{E}(W_2) - \mathcal{E}(W_1) \leq \delta;$
\end{enumerate}
then $C \geq F(\lfloor \sqrt{\delta^{-1}}/10p \rfloor-1)$ where $F: \mathbb{N} \rightarrow \mathbb{N}$ is a function defined by $F(0) = 0$ and $F(i+1) = \twr\left(\lfloor\sqrt{\delta}/8p\epsilon(F(i))\rfloor\right)$.
\end{theorem}

\begin{proof}
    Let $f$ be defined as at the end of Section \ref{section:construction}. Writing $t = \lfloor \sqrt{\delta^{-1}}/10p \rfloor$, suppose that there is an integer $1 \leq k \leq t-1$ such that $C \leq D_{\phi(k, h_k-9p)}$. Take $k$ to be minimal with this property, so that, in particular, $C \geq D_{\phi(k-1, h_{k-1}-9p)}$ (or $C \geq D_0$ if $k=1$). Then $\epsilon(C) \leq \epsilon(D_{\phi(k-1, h_{k-1}-9p)}) = \epsilon_{k}$ and therefore $\mathcal{P}(W_2)$ is an $\epsilon_{k}$-regular partition for $f$.
    
    Since $w_{k, h_k} \geq 8p \epsilon_k$ by construction, Proposition \ref{prop:large-bias} implies that $W_2 \leqslant H_{\phi(k, h_k)}$. Moreover, $\mathcal{E}(W_1) <\mathcal{E}(H_{\phi(k, h_k-9p)}) + 8 w_{k, h_k-9p}^2$ by Proposition \ref{prop:energy-start} so $$\mathcal{E}(W_2) - \mathcal{E}(W_1) > \mathcal{E}(H_{\phi(k, h_k)}) - \mathcal{E}(H_{\phi(k, h_k-9p)}) - 8 w_{k, h_k-9p}^2.$$
    However, Proposition \ref{prop:energy-middle} gives
        $$\mathcal{E}(H_{\phi(k, h_k)}) - \mathcal{E}(H_{\phi(k, h_k-9p)}) \geq \frac1p w_{k, h_k-1}^2 + \ldots + \frac1p w_{k, h_k-9p}^2 = 9w_{k, h_k-9p}^2,$$
    and therefore $\mathcal{E}(W_2) - \mathcal{E}(W_1) > w_{k, h_k-9p}^2 \geq \delta$, which is a contradiction. As a result, we must conclude that such a $k$ does not exist, i.e.~$C \geq D_{\phi(t-1, h_{t-1}-9p)}$.
    
    It remains to show that $D_{\phi(t-1, h_{t-1}-9p)} \geq F(t-1)$. Recall that $h_1 \geq 10p$ and that $D_k \geq \twr(k)$ for any $k \geq 1$. Then for any $i \geq 2$,
    \vspace{-0.3cm}
    $$D_{\phi(i, h_i-9p)} \geq \twr\left(h_1+\sum_{j=2}^{i} h_j-9p\right) \geq \twr\left(h_i\right).  \vspace{-0.3cm}$$
    Suppose that $D_{\phi(i, h_i-9p)} \geq F(i)$ for some $1 \leq i \leq t-1$. Then $h_{i+1} = \lfloor\sqrt{\delta}/8p\epsilon_{i+1}\rfloor$ is greater than $\lfloor\sqrt{\delta}/8p\epsilon(F(i))\rfloor$, so
    $$D_{\phi(i+1, h_{i+1}-9p)} \geq \twr(h_{i+1}) \geq \twr\left(\lfloor\sqrt{\delta}/8p \epsilon(F(i))\rfloor\right) = F(i+1).$$
    As a consequence, $D_{\phi(t-1, h_{t-1}-9p)} \geq F(t-1)$ holds by induction, with $D_{\phi(1, h_1-9p)} \geq F(1)$ taken as the base case.
\end{proof}

\noindent As the definition of $F$ involves iterating a tower function, $F$ indeed has wowzer-type growth in its parameters. In fact, the following restatement of Theorem \ref{thm:sarl-intro} from the introduction is an immediate consequence of Theorem \ref{thm:sarl-lowerbound}.

\begin{mainA}[Wowzer-type lower bound on $C_{sarl}^\bracketed{1}$]
    Fix $0 < \delta \leq 1/20p$ and $\epsilon: \mathbb{N} \rightarrow (0,1)$ such that $\epsilon(d) \leq \sqrt{\delta}/(80p^2(d+1))$. There exists a function $f: \vecsp{p}{n} \rightarrow [0,1]$ for which the following holds. If $W_2 \leqslant W_1 \leqslant \vecsp{p}{n}$ are subspaces of codimensions $C_2$ and $C_1$ respectively satisfying
    \begin{enumerate}[label=(\roman*)]
        \item $\mathcal{P}(W_2)$ is $\epsilon(C_1)$-regular for $f$;
        \item $\mathcal{E}(W_2) - \mathcal{E}(W_1) \leq \delta;$
    \end{enumerate}
then $C_2, C_1 > \wwz(\lfloor \sqrt{\delta^{-1}}/10p \rfloor)$. In particular, $C_{sarl}^\bracketed{1}(p, \delta, \epsilon) > \wwz(\lfloor \sqrt{\delta^{-1}}/10p \rfloor)$.
\end{mainA}
\begin{proof}
    By Theorem \ref{thm:sarl-lowerbound} applied with $\delta$ and $\epsilon$, $C_1 \geq F(t - 1)$ where $t = \lfloor \sqrt{\delta^{-1}}/10p \rfloor$. Firstly, note that $F(1) \geq \twr(10p) \geq \twr(2) = \wwz(2)$. On the other hand, if $F(i) > \wwz(i+1)$ for some $1 \leq i < t$, then
    \vspace{-0.1cm}
    $$F(i+1) \geq \twr\left(10p F(i)\right) \geq  \twr\left(\wwz(i+1)\right) = \wwz(i+2).$$
    Hence $C_2 \geq C_1 \geq F(t - 1) > \wwz(t)$ by induction.
\end{proof}

\section{Size of the linear layer in higher-order regularity partitions}
\label{sec:linear-layer}

As noted in the introduction, Fourier uniformity allows us to count the number of solutions to linear systems of true complexity 1. The concept of true complexity of a linear system was developed by Gowers and Wolf \cite{true-complexity}. Informally, it is defined as the smallest integer $s$ (or $\infty$) such that the number of solutions to the given linear system in any set $A \subseteq \vecsp{p}{n}$ (or, more generally, under a function $f:\vecsp{p}{n} \rightarrow [0,1]$) is `controlled' by the Gowers uniformity norm $\norm{\id{A}}_{U^{s+1}}$ \cite{gowers-norm}. The latter will not be defined here, but it is a well-known fact that $\norm{f}_{U^2} = \norm{\widehat{f}}_4$ (\cite[Lemma 2.4]{gowers-decompositions}), which explains why Fourier uniformity controls systems of complexity 1.

Indeed, there are higher-order analogues of the arithmetic regularity lemma (Theorem \ref{theorem:arl}) corresponding to the norm $U^{s+1}$ for each $s \geq 2$ (see, for instance, \cite{montreal}, \cite[Section 4]{higher-fourier}, \cite{decomposition-lemmas} or \cite{gowers-wolf-higher-degree}). Like Theorem \ref{theorem:arl}, these higher-order arithmetic regularity lemmas provide a partition of $\vecsp{p}{n}$ that is `regular' in some sense for the given function $f$. The goal of this section is to demonstrate that the `linear layer' of such partitions must still be of tower-type size in some cases, despite the additional features arising in this setting; in fact, this matches the upper bound in the order of growth arising from the proof of the quadratic regularity lemma (see Appendix \ref{appendix:upper-bounds}). For ease of exposition, all proofs in this section are presented for the quadratic setting $s=2$, but arithmetic regularity lemmas of orders higher than 2 may be treated similarly.

Where Theorem \ref{theorem:arl} produces a partition of $\vecsp{p}{n}$ into cosets of a subspace, which may be viewed as a linearly-structured partition, a quadratic arithmetic regularity lemma utilises a quadratically-structured partition defined as follows.

\begin{definition}[Quadratic factor]
    \label{def:quadratic-factor}
    Given polynomials $P_1, \ldots, P_D:\vecsp{p}{n} \rightarrow \vecsp{p}{}$ of degree at most $2$, a \emph{quadratic factor $\mathcal{B}$} of complexity $D$ is a partition of $\vecsp{p}{n}$ into the simultaneous level sets of $(P_1, \ldots, P_D)$, referred to as \emph{atoms}. As such, an atom $B$ of $\mathcal{B}$ has the form
    $$B = \set{x \in\vecsp{p}{n}\st (P_1(x), \ldots, P_D(x))=\vecline{c}},$$
    where $\vecline{c} \in \vecsp{p}{D}$ is the \emph{label of $B$}. The set of all atoms of $\mathcal{B}$ is denoted by $\At(\mathcal{B})$.
    
    Additionally, we will write $\mathcal{B}[1]$ for the \emph{linear layer of $\mathcal{B}$}, which is the coarser factor defined by the linear polynomials $(P_i: \deg(P_i) = 1)$, and $\mathcal{B}[2]$ for the \emph{quadratic layer} defined by the quadratic polynomials $(P_i: \deg(P_i) = 2)$.
\end{definition}
\noindent \textbf{Note:} If $L_1, \ldots, L_{\ell}$ are the polynomials defining $\mathcal{B}[1]$, then $\mathcal{B}[1]$ is a partition of $\vecsp{p}{n}$ into cosets of the subspace $H = \set{x \in \vecsp{p}{n}: L_i(x)=0}$. In this way, quadratic factors may be seen as partitions into cosets that are further refined into quadratically structured parts.

In applications, it is convenient to work with quadratic factors of \emph{high rank} as this ensures that all atoms have approximately the same size. If $\mathcal{B}$ is a quadratic factor and $\mathcal{B}[2]$ is defined by quadratic polynomials $Q_1, \ldots, Q_q$, then the rank of $\mathcal{B}$ is the minimum rank of any non-zero linear combination of $Q_1, \ldots, Q_q$, i.e.
$$\rank(\mathcal{B}) = \min_{\lambda \in \vecsp{p}{q}\backslash\set{0}}\rank(\sum_{i=1}^q \lambda_i Q_i ),$$
where the rank of a quadratic polynomial $Q$ is simply the matrix rank of the unique $n \times n$ matrix $M$ over $\vecsp{p}{n}$ such that $Q(x) = x^T M x + L(x)$ for some linear polynomial $L$. The following result is a combination of \cite[Lemma 3.1]{montreal} and \cite[Lemma 4.2]{montreal}.

\begin{lemma}[High rank implies equidistribution \cite{montreal}]
\label{lemma:atom-size}
     Let $\mathcal{B}$ be a quadratic polynomial of complexity $D$ and rank $r$. If $B$ is an atom of $\mathcal{B}$, then
     \begin{enumerate}[label = (\roman*)]
         \item \label{prpt:quadratic-gaussian} for any linear polynomial $L$, $\babs{\expct_{x \in B} e(L(x))} \leq p^{-r/2}$;

         \item \label{prpt:atom-size} if $r \geq 2(D+1)$, then $p^{-D}/2 \leq \abs{B}/\abs{\vecsp{p}{n}} \leq 3 p^{-D}/2$.
     \end{enumerate}
\end{lemma}

With these concepts defined, it is now possible to describe quadratic regularity partitions precisely. In the definition below, $\expct(f|\mathcal{B})$ denotes the projection of $f$ onto $\mathcal{B}$ where $\expct(f|\mathcal{B})(x)$ is equal to the average of $f$ on the atom of $\mathcal{B}$ containing $x$.
\begin{definition}[Quadratic regularity partition]
\label{def:quadratic-regularity-partition}
Fix $\delta > 0$, two non-decreasing functions $\omega, R: \mathbb{N} \rightarrow \mathbb{N}$, and let $f:\vecsp{p}{n} \rightarrow [0, 1]$ be a function. A quadratic factor $\mathcal{B}$ of complexity $D$ is a \emph{$(\delta, \omega, R)$-quadratic regularity partition for $f$} if $\rank(\mathcal{B}) \geq R(D)$ and there is a function $f_{err}:\vecsp{p}{n} \rightarrow [-1,1]$ such that $\norm{f_{err}}_{L_2} < \delta$ and $\norm{f- \expct(f|\mathcal{B}) - f_{err}}_{U^3} < 1/\omega(D)$.
\end{definition}
\noindent To the reader unfamiliar with higher-order arithmetic regularity lemmas, Definition \ref{def:quadratic-regularity-partition} may seem to not bear much resemblance to the regular partition defined in Definition \ref{def:regularity}. However, it is possible to express the latter in a form much like Definition \ref{def:quadratic-regularity-partition}: the $U^2$-norm would replace the $U^3$-norm, $\mathcal{B}$ would be a `linear' factor, and there is no need for a rank function; a `linear regularity partition' thus defined may be translated to the usual form in terms of Fourier-uniformity (see \cite[Lemma 2.10]{tao-structure}).

It turns out that that a similar translation can be carried out in the quadratic case. Just as there is a connection between the $U^2$-norm and the Fourier transforms of $f$, the $U^3$-norm in some way corresponds to bias with respect to quadratic polynomial phases (see, for instance, the inverse theorem for the $U^3$-norm \cite[Lecture 2]{montreal}). This fact informs the definition below, with the notation for quadratic bias taken from \cite[Definition 2.1]{inverse-u3}.

\begin{definition}[$\epsilon$-quadratically unbiased on $B$]
\label{def:unbias}
Let $\mathcal{B}$ be a quadratic factor and let $\polys{2}$ denote the set of polynomials of degree at most $2$ in $\vecsp{p}{n}$. Given an atom $B$ of $\mathcal{B}$ and a function $f: \vecsp{p}{n} \rightarrow [0,1]$, the \emph{quadratic bias of $f$ on $B$} is defined as
$$\norm{f}_{u^{3}(B)} = \sup_{P \in \polys{2}} \babs{\expct_{x \in B} f(x) e(P(x))}.$$
For $\epsilon > 0$, $f$ is said to be \emph{$\epsilon$-quadratically unbiased on $B$} if $\norm{f-\alpha_B}_{u^{3}(B)} \leq \epsilon$, where $\alpha_B$ is the density of $f$ on $B$.
\end{definition}

Quadratic bias is a natural generalisation of Fourier uniformity, which is itself a measure of linear bias. It is now possible to bring Definition \ref{def:quadratic-regularity-partition} more in line with the regular partitions defined in Definition \ref{def:regularity} via the following lemma, whose proof may be found in Appendix \ref{appendix:proof-of-lemma}. Note that $\omega(d) \geq \delta^{-2/3} p^{d}$ is not an unreasonable assumption, since applications typically require $\omega(d) = 2^{O(d)}$ (for example, see the proof of Theorem 4.1 in \cite{montreal} or the choice of parameters in the proof of Theorem 5.10 in \cite{trans-invariant}).

\begin{lemma}[Unbiased quadratic regularity partition]
\label{lemma:higher-arl-atom-version}
  Fix $\delta > 0$, two non-decreasing functions $\omega, R: \mathbb{N} \rightarrow \mathbb{N}$ such that $\omega(d) \geq \delta^{-2/3} p^{d}$ and $R(d) \geq 2(d+1)$, and let $f:\vecsp{p}{n} \rightarrow [0, 1]$ be a function. If $\mathcal{B}$ is a $(\delta, \omega, R)$-quadratic regularity partition, then for all but a $2 \delta^{2/3}$-proportion of atoms $B$ of $\mathcal{B}$, $f$ is $3\delta^{2/3}$-quadratically unbiased on $B$.
\end{lemma}

\noindent We will now show that if $\mathcal{B}$ is a quadratic regularity partition for a function $f:\vecsp{p}{n} \rightarrow [0,1]$, then $f$ is Fourier-uniform on almost all cosets in the underlying partition $\mathcal{B}[1]$, i.e.~$\mathcal{B}[1]$ is a regular partition for $f$ in the sense of Definition \ref{def:regularity}. The main result of the section then follows easily as a corollary.

\begin{prop}
    \label{prop:B1-is-regular}
     Fix $\delta > 0$, two non-decreasing functions $\omega, R: \mathbb{N} \rightarrow \mathbb{N}$ such that $\omega(d) \geq \delta^{-2/3} p^{d}$ and $R(d) \geq 2(d+1 + \log_p(\delta^{1/3}))$, and let $f:\vecsp{p}{n} \rightarrow [0, 1]$ be a function. If $\mathcal{B}$ is a $(\delta, \omega, R)$-quadratic regularity partition, then $\mathcal{B}[1]$ is $7\delta^{1/3}$-regular for $f$ in the sense of Definition \ref{def:regularity}.
\end{prop}
\begin{proof}
    Let $D$ denote the complexity of $\mathcal{B}$, and let $H$ be the subspace such that $\mathcal{B}[1]$ is a partition into cosets of $H$. Additionally, for each $c \in \vecsp{p}{n}$, write $\At_c$ for the set of atoms of $\mathcal{B}$ contained in the coset $H+c$, noting that $\At_c$ has the same size for every $c$. By Lemma \ref{lemma:higher-arl-atom-version}, for all but a $2 \delta^{2/3}$-proportion of atoms $B$ of $\mathcal{B}$, $f$ is $3\delta^{2/3}$-quadratically unbiased on $B$. By averaging, this implies that for all but a $2 \delta^{1/3}$-proportion of $c \in \vecsp{p}{n}$, $f$ is $3\delta^{2/3}$-quadratically unbiased on all but a $\delta^{1/3}$-proportion of atoms in $\At_c$.

    \begin{claim}
        Let $c \in \vecsp{p}{n}$ be such that $f$ is $3\delta^{2/3}$-quadratically unbiased on all but a $\delta^{1/3}$-proportion of atoms in $\At_c$. Then $f$ is $7\delta^{1/3}$-Fourier uniform on $H+c$.
    \end{claim}
    \begin{proof}[Proof of Claim]
        \renewcommand{\qedsymbol}{}
        Let $\alpha_{c}$ denote the density of $f$ on $H+c$ and, given an atom $B \in \At_c$, let $\beta_B$ denote the density of $f$ on $B$.  Fix any $r \in \vecsp{p}{n}$ and define the linear polynomial $L_r(x) = r^T x$.
        Writing $F$ for the balanced function $f-\alpha_c$ on $H+c$, the Fourier transform of $f$ on $H+c$ at $r$ can be rewritten as
        $$\widehat{F}(r) = \expct_{x \in H+c} F(x) e_p(L_r(x)) = \sum_{B \in \At_c} \frac{\abs{B}}{\abs{H}} \expct_{x \in B} F(x) e_p(L_r(x)).$$
        Writing $F_B = f - \beta_B$ for the balanced function of $f$ on $B$, this becomes
        \begin{equation}
            \label{eq:F_Bs-and-expct}
            \widehat{F}(r) = \sum_{B \in \At_c} \frac{\abs{B}}{\abs{H}} \expct_{x \in B} F_B(x) e_p(L_r(x)) + \sum_{B \in \At_c} \beta_B \frac{\abs{B}}{\abs{H}} \expct_{x \in B} e_p(L_r(x)).
        \end{equation}
        The second term here can be bounded by the triangle inequality and Lemma \ref{lemma:atom-size}\ref{prpt:quadratic-gaussian} so 
        \begin{equation}
            \label{eq:sum-of-gaussian-expct}
            \left| \sum_{B \in \At_c} \beta_B \frac{\abs{B}}{\abs{H}}  \expct_{x \in B} e_p(L_r(x))  \right| \leq \sum_{B \in \At_c} \beta_B \frac{\abs{B}}{\abs{H}} \Big|\expct_{x \in B} e_p(L_r(x))\Big| \leq p^{-R(D)/2} \leq \delta^{1/3},
        \end{equation}
        where the last inequality uses the fact that $\sum_{B \in \At_c} \beta_B {\abs{B}}/{\abs{H}} = \alpha_c \leq 1$.
        
        For the first term of equation \eqref{eq:F_Bs-and-expct}, observe that for all atoms $B$ on which $f$ is $3\delta^{2/3}$-quadratically unbiased, $\babs{\expct_{x \in B} F_B(x) e_p(L_r(x))} \leq 3\delta^{2/3}$, and there is only a $\delta^{1/3}$-proportion of atoms in $\At_c$ for which this may not be the case. Moreover, for all $B \in \At_c$, $\abs{B}/\abs{H} \leq 3 \abs{\At_c}^{-1}/2$ by Lemma \ref{lemma:atom-size}\ref{prpt:atom-size}. As a result,
        \begin{equation}
            \label{eq:sum-of-unbiased-expct}
            \left| \sum_{B \in \At_c} \frac{\abs{B}}{\abs{H}} \expct_{x \in B} F_B(x) e_p(L_r(x))  \right| \leq \frac32 \expct_{B \in \At_c} \Big| \expct_{x \in B} F_B(x) e_p(L_r(x)) \Big| \leq \frac32 \Big(3\delta^{2/3} + \delta^{1/3}\Big).
        \end{equation}

        Finally, apply the triangle inequality to equation \eqref{eq:F_Bs-and-expct} and use the bounds from \eqref{eq:sum-of-gaussian-expct} and \eqref{eq:sum-of-unbiased-expct} to deduce that $\babs{\widehat{F}(r)} \leq 6\delta^{1/3} + \delta^{1/3} \leq 7\delta^{1/3}$, which proves the claim.
    \end{proof}

    \noindent As already established, the claim applies to all but a $2\delta^{1/3}$-proportion of $c \in \vecsp{p}{n}$, which is certainly less than a $7\delta^{1/3}$-proportion. In particular, it follows that $\mathcal{B}[1]$ is a $7\delta^{1/3}$-regular partition, as required.
\end{proof}

\noindent Combining Proposition \ref{prop:B1-is-regular} with Corollary \ref{corr:hosseini-for-p} immediately gives Theorem \ref{thm:main-linear-layer}.

\begin{mainB}[Tower-type lower bound on the linear layer]
    Fix $\delta > 0$, two non-decreasing functions $\omega, R: \mathbb{N} \rightarrow \mathbb{N}$ such that $\omega(d) \geq \delta^{-2/3} p^{d}$ and $R(d) \geq 2(d+1 + \log_p(\delta^{1/3}))$. There exists a function $f:\vecsp{p}{n} \rightarrow [0, 1]$ such that, if $\mathcal{B}$ is a $(\delta, \omega, R)$-quadratic regularity partition for $f$, then $\mathcal{B}[1]$ is a partition into cosets of a subspace of codimension at least $\twr(\lfloor\delta^{-1/3}/60p\rfloor - 1)$.
\end{mainB}
\noindent Note that a similar result may be obtained for arithmetic regularity lemmas of orders higher than $2$ by following the same proof, with only small technical modifications.

\printbibliography

\appendix
\section{Proof of Lemma \ref{lemma:higher-arl-atom-version}}
\label{appendix:proof-of-lemma}

Recall the definition of quadratic bias $\|\cdot \|_{u^{3}(B)}$ in Definition \ref{def:unbias}. It is a well-known fact that quadratic bias on the whole of $\vecsp{p}{n}$ is controlled by the Gowers $U^3$-norm \cite[Equation (2.2)]{inverse-u3}: specifically, for any $f: \vecsp{p}{n} \rightarrow \mathbb{C}$,
\begin{equation}
    \label{eq:bias-gowers-bound}
    \norm{f}_{u^{3}(\vecsp{p}{n})} \leq \norm{f}_{U^{3}}.
\end{equation}
The following lemma leverages this fact to turn global uniformity in terms of the $U^3$-norm into a lack of local bias in the sense of Definition \ref{def:unbias}.

\begin{lemma}[Uniformity implies lack of local bias]
\label{lemma:uniform-to-bias}
Fix $\delta, \eta > 0$. Let $\mathcal{B}$ be a quadratic factor of complexity $D$, and let $f: \vecsp{p}{n} \rightarrow [0,1]$ be a function. Suppose that $B$ is an atom of $\mathcal{B}$ such that
there is a function $f_{err}: \vecsp{p}{n} \rightarrow [-1,1]$ satisfying
\begin{itemize}
    \item $\norm{f - \expct_{}[f|\mathcal{B}] - f_{err}}_{U^3} \leq \eta$, and
    \item $\expct_{x \in B} \abs{f_{err}(x)}^2 \leq \delta^2$.
\end{itemize}
Then $f$ is $\epsilon$-unbiased on $B$ for some $\epsilon \leq \delta + \eta {\abs{\vecsp{p}{n}}}/{\abs{B}}.$
\end{lemma}
\begin{proof}
Let $P$ be a polynomial of degree at most $2$. With $F_B$ denoting the balanced function $f-\expct_{x \in B} f(x)$ on $B$, the triangle inequality gives
\begin{align}
\label{eq:globalunf-to-localunbias}
\babs{\expct_{x \in B} F_B(x) e(P(x))} &\leq \babs{\expct_{x \in B} f_{err}(x) e(P(x))} + \babs{\expct_{x \in B} (F_B - f_{err})(x) e(P(x))}.
\end{align}
The first term here can be bounded by the square root of $\expct_{x \in B} \abs{f_{err}(x)}^2$ via the Cauchy-Schwarz inequality, and the second term can be handled with the following claim.
\begin{claim}
    For any function $F:\vecsp{p}{n} \rightarrow [-1,1]$,
    $\babs{\expct_{x \in B} F(x) e(P(x))} \leq \norm{F}_{U^3} {\abs{\vecsp{p}{n}}}/{\abs{B}}.$
\end{claim}
\begin{proof}[Proof of Claim]
     Let $\underline{a} = (a_1, \ldots, a_D) \in \vecsp{p}{D}$ denote the label of $B$. Then $\id{B}(x)$ may be rewritten as the product of $\id{P_i(x)=a_i}(x) = \expct_{\lambda_i \in \vecsp{p}{}} e(\lambda_i (P_i(x) - a_i))$, resulting in
    \vspace{-0.2cm}
    \begin{align*}
        \expct_{x \in B} F(x) e(P(x)) &= \frac{\abs{\vecsp{p}{n}}}{\abs{B}} \expct_{x \in \vecsp{p}{n}} \Big{[} F(x) e(P(x)) \prod_{i=1}^D \expct_{\lambda_i \in \vecsp{p}{}} e(\lambda_i (P_i(x) - a_i))\Big{]}\\
        &= \frac{\abs{\vecsp{p}{n}}}{\abs{B}} \expct_{\underline{\lambda} \in \vecsp{p}{D}} \expct_{x \in \vecsp{p}{n}} F(x) e\Big{(}P(x) + \sum_{i=1}^D \lambda_i (P_i(x) - a_i)\Big{)}.
    \end{align*}
    The claim follows by applying \eqref{eq:bias-gowers-bound}.
\end{proof}

\noindent As a consequence, $\babs{\expct_{x \in B} (F_B - f_{err})(x) e(P(x))} \leq \norm{F_B - f_{err}}_{U^3}{\abs{\vecsp{p}{n}}}/{\abs{B}}$. Substituting this into equation \eqref{eq:globalunf-to-localunbias} and recalling that, by assumption, $\norm{F_B - f_{err}}_{U^{s+1}} \leq \eta$ and $\expct_{x \in B} \abs{f_{err}(x)}^2 \leq \delta^2$ completes the proof.
\end{proof}

The proof of Lemma \ref{lemma:higher-arl-atom-version} now proceeds by showing that there are many atoms $B$ on which $\expct_{x \in B} \abs{f_{err}(x)}^2 \leq \delta^2$ and using Lemma \ref{lemma:uniform-to-bias}.

\begin{replemma}{lemma:higher-arl-atom-version}[Unbiased quadratic regularity partition]
 Fix $\delta > 0$, two non-decreasing functions $\omega, R: \mathbb{N} \rightarrow \mathbb{N}$ such that $\omega(d) \geq \delta^{-2/3} p^{d}$ and $R(d) \geq 2(d+1)$, and let $f:\vecsp{p}{n} \rightarrow [0, 1]$ be a function. If $\mathcal{B}$ is a $(\delta, \omega, R)$-quadratic regularity partition, then for all but a $2 \delta^{2/3}$-proportion of atoms $B$ of $\mathcal{B}$, $f$ is $(3\delta^{2/3}, 2)$-unbiased on $B$.
\end{replemma}
\begin{proof}
By Definition \ref{def:quadratic-regularity-partition}, $\rank(\mathcal{B}) \geq 2(D+1)$ and there is a function $f_{err}:\vecsp{p}{n} \rightarrow [-1,1]$ such that $\norm{f_{err}}_{L_2} < \delta$ and $\norm{f - \expct_{}[f|\mathcal{B}] - f_{err}}_{U^{s+1}} < 1/\omega(D)$.

\begin{nclaim}
    For all but a $2 \delta^{2/3}$-proportion of atoms $B$ of $\mathcal{B}$, $\expct_{x \in B}\abs{f_{err}(x)}^2 \leq \delta^{4/3}$.
\end{nclaim}
\begin{proof}[Proof of Claim.]
    \renewcommand{\qedsymbol}{}
    Rewrite $\norm{f_{err}}_{L_2}^2$ as a sum over the atoms of $\mathcal{B}$, i.e.
    \begin{equation}
        \label{eq:arl-l2-over-atoms}
        \norm{f_{err}}_{L_2}^2 = \sum_{B \in \At(\mathcal{B})} \expct_{x \in \vecsp{p}{n}} \abs{f_{err}(x)}^2 \id{B}(x) = \sum_{B \in \At(\mathcal{B})} \frac{\abs{B}}{\abs{\vecsp{p}{n}}} \expct_{x \in B}\abs{f_{err}(x)}^2.
    \end{equation}
    By Lemma \ref{lemma:atom-size}\ref{prpt:atom-size}, $\abs{B}/\abs{\vecsp{p}{n}} \geq p^{-D}/2$, and $\norm{f_{err}}_{L_2}^2 \leq \delta^2$ by assumption, so
    $$\expct_{B \in \At(\mathcal{B})} \bigg[\expct_{x \in B}\abs{f_{err}(x)}^2\bigg] \leq 2 \norm{f_{err}}_{L_2}^2 < 2\delta^2.$$
    By averaging, $\expct_{x \in B}\abs{f_{err}(x)}^2 \leq \delta^{4/3}$ for all but a $2\delta^{2/3}$-proportion of $B$, as required.
\end{proof}

\noindent Now observe that for each $B$ satisfying the conclusion of the claim, Lemma \ref{lemma:uniform-to-bias} implies that $f$ is $(\epsilon, 2)$-unbiased on $B$ with $\epsilon \leq \delta^{2/3} + \omega(D)^{-1}{\abs{\vecsp{p}{n}}}/{\abs{B}}.$ Since $\omega(D) \geq \delta^{-2/3}p^{D}$ by assumption and $\abs{B}/\abs{\vecsp{p}{n}} \geq p^{-D}/2$ by Lemma \ref{lemma:atom-size}\ref{prpt:atom-size}, $\epsilon \leq \delta^{2/3} + 2 \delta^{2/3} {p^{-D}} p^D \leq 3 \delta^{2/3}$, which gives the required conclusion.
\end{proof}

\section{An upper bound for quadratic regularity partitions}
\label{appendix:upper-bounds}

Section \ref{sec:linear-layer} of this paper concerns quadratic regularity partitions (recall Definition \ref{def:quadratic-regularity-partition}), showing that for some functions $f:\vecsp{p}{n} \rightarrow [0,1]$, any sufficiently high-rank $(\delta, \omega, R)$-quadratic regularity partition must have a linear layer of size at least tower-type in $\delta^{-O(1)}$. On the other hand, for any function $f$ and any choice of parameters, one can always find a $(\delta, \omega, R)$-quadratic regularity partition for $f$ by the quadratic arithmetic regularity lemma (e.g.~\cite[Proposition 3.12]{montreal}).

\begin{theorem}[Quadratic arithmetic regularity lemma]
\label{thm:higher-arl}
Fix $\delta > 0$ and two non-decreasing functions $\omega, R: \mathbb{N} \rightarrow \mathbb{N}$. There is a constant $C_{qarl} = C_{qarl}(\delta, \omega, R)$ such that the following holds. For all functions $f:\vecsp{p}{n} \rightarrow [0, 1]$, there exists a $(\delta, \omega, R)$-quadratic regularity partition $\mathcal{B}$ of complexity $D \leq C_{qarl}$.
\end{theorem}

\noindent It is not hard to see, by following a standard proof of Theorem \ref{thm:higher-arl}, that the upper bound on the size of such a partition is tower-type in $\delta^{-O(1)}$, with details provided in this appendix for completeness.

Proofs of regularity lemmas are typically energy increment arguments (see, for instance, \cite{tower-type} or \cite{montreal}), with energy defined as in Definition \ref{def:energy}. In the quadratic setting, the proof proceeds as follows. Start with a trivial partition $\mathcal{B}_0$, noting that it has infinite rank, and write $\eta_t = 1/\omega(t)$ for each $t \in \mathbb{N}$. If $\mathcal{B}_0$ satisfies $\norm{f- \expct(f|\mathcal{B}_0)}_{U^{3}} < \eta_0$, then we are done by taking $f_{err} = 0$; otherwise, since the $U^{3}$-norm is large, $f- \expct(f|\mathcal{B}_0)$ must have relatively large bias with a polynomial of degree at most 2 by the inverse theorem for the $U^3$-norm \cite{tao-ziegler}. In fact, the dependence between the size of the $U^3$-norm and the quadratic bias is polynomial as a consequence of the recent work of Gowers, Green, Manners, and Tao \cite{pfr} on the polynomial Freiman-Rusza conjecture (see \cite[Corollary 1.6]{pfr} as well as \cite{pfr-implies-inverse-theorem}).

\begin{theorem}[Inverse theorem for the $U^3$-norm]
\label{theorem:gowers-inverse}
For all $\eta > 0$ and all $g: \vecsp{p}{n} \rightarrow [-1,1]$ the following holds. If $\norm{g}_{U^3} > \eta$, then there is a polynomial $P:\vecsp{p}{n} \rightarrow \vecsp{p}{}$ of degree at most $2$ such that
$\babs{\expct_{x \in \vecsp{p}{n}}g(x) e(P(x))} > \eta^{O(1)}.$
\end{theorem}

Apply Theorem \ref{theorem:gowers-inverse} to $f - \expct(f|\mathcal{B}_0)$ to obtain a polynomial $P$ of degree at most 2, and let $\mathcal{B}'_0$ denote the quadratic factor defined by $P$. It can be shown \cite[Lemma 3.8]{montreal} that this step increases the energy of the underlying partition by $\eta_0^{O(1)}$.

\begin{lemma}[Energy increment \cite{montreal}]
Fix $\eta > 0$ and let $\mathcal{B}$ be a quadratic factor of complexity $D$. If $g:\vecsp{p}{n} \rightarrow [-1,1]$ is a function such that $\norm{g - \expct(g|\mathcal{B})}_{U^3} \geq \eta,$
then there is a quadratic factor $\mathcal{B}'$ refining $\mathcal{B}$ such that $\mathcal{B}'$ has complexity $D+1$ and
$\mathcal{E}(\mathcal{B}') \geq \mathcal{E}(\mathcal{B}) + \eta^{O(1)}.$
\end{lemma}

If it is also the case that $\norm{f- \expct(f|\mathcal{B}'_0)}_{U^3} > \eta_0$, then $f- \expct(f|\mathcal{B}'_0)$ must correlate with another quadratic polynomial by Theorem \ref{theorem:gowers-inverse}. Add this polynomial into the definition of $\mathcal{B}'_0$ to obtain a new factor $\mathcal{B}''_0$ such that $\mathcal{E}(\mathcal{B}''_0) \geq \mathcal{E}(\mathcal{B}'_0) + \eta_0^{O(1)}$, and repeat again with $\mathcal{B}''_0$. Since energy can only take values between 0 and 1 by Lemma \ref{lemma:energy-props}\ref{energyprop:bounded}, such an iteration must terminate in at most $\eta_0^{-O(1)}$ steps, resulting in a quadratic factor $\mathcal{I}_{0}$ of complexity at most $\eta_0^{-O(1)}$ such that $\norm{f- \expct(f|\mathcal{I}_0)}_{U^3} \leq \eta_0$.

In fact, if we take $f_{err}^\bracketed{0} = \expct(f|\mathcal{I}_0) - \expct(f|\mathcal{B}_0)$, then $\norm{f- \expct(f|\mathcal{B}_0) - f_{err}^\bracketed{0}}_{U^3} \leq \eta_0$. Since $\mathcal{B}_0$ has infinite rank, this would make it a $(\delta, \omega, R)$-quadratic regularity partition for $f$ \textit{provided} that $\norm{f_{err}^\bracketed{0}}_{L_2}^2 \leq \delta^2$. In the case that the latter does not hold, observe that $\norm{f_{err}^\bracketed{0}}_{L_2}^2 = \mathcal{E}(\mathcal{I}_0) - \mathcal{E}(\mathcal{B}_0)$ by Lemma \ref{lemma:energy-props}\ref{energyprop:pythagoras}, so $\norm{f_{err}^\bracketed{0}}_{L_2}^2  > \delta^2$ simply gives us another energy increment. The strategy, then, is to repeat the whole argument with $\mathcal{I}_0$ in place of $\mathcal{B}_0$, and so on, establishing an energy increment $\delta^2$ at every step until, in at most $\delta^{-2}$ steps, we arrive at a $(\delta, \omega, R)$-quadratic regularity partition.

The only caveat is that $\mathcal{I}_0$ may not have sufficiently high rank anymore. This can be resolved by an additional refining step, such as \cite[Lemma 3.11]{montreal}, restated below.
\begin{lemma}[Making quadratic factors high-rank \cite{montreal}]
    \label{lemma:high-rank-refinement-qdr}
    Let $R: \mathbb{N} \rightarrow \mathbb{N}$ be a non-decreasing function. There exists a function $\phi_{R}: \mathbb{N} \rightarrow \mathbb{N}$ satisfying the following. For every quadratic factor $\mathcal{B}$ of complexity $D$, there is a quadratic factor $\mathcal{B}'$ refining $\mathcal{B}$ such that $\mathcal{B}'$ has complexity $D' \leq \phi_{R}(D)$ and $\rank(\mathcal{B}') \geq R(D')$.
\end{lemma}

Let $\mathcal{B}_1$ be the result of applying Lemma \ref{lemma:high-rank-refinement-qdr} to $\mathcal{I}_0$, so that $\mathcal{B}_1$ has rank at least $R(D_1)$ and complexity $D_1 \leq \phi_{R}(\eta_0^{-O(1)})$. Crucially, by Lemma \ref{lemma:energy-props}\ref{energyprop:energy-increase}, there is still an energy increment of $\mathcal{E}(\mathcal{B}_1) - \mathcal{E}(\mathcal{B}_0) \geq \delta^2$.

Now we may truly carry out the proposed strategy. Starting with $\mathcal{B}_1$ in place of $\mathcal{B}_0$, argue as before to obtain a quadratic factor $\mathcal{I}_1$ refining $\mathcal{B}_1$ such that, firstly, $\mathcal{I}_1$ has complexity at most $D_1+\eta_{D_1}^{-O(1)}$ and, secondly, $\norm{f- \expct(f|\mathcal{B}_1) - f_{err}^\bracketed{1}}_{U^3} \leq \eta_{D_1}$, where $f_{err}^\bracketed{1} =\expct(f|\mathcal{I}_1) - \expct(f|\mathcal{B}_1)$. If $\norm{f_{err}^\bracketed{1}}_{L_2}^2  > \delta^2$, apply Lemma \ref{lemma:high-rank-refinement-qdr} to $\mathcal{I}_1$ to obtain a quadratic factor $\mathcal{B}_2$ of complexity $D_2 \leq \phi_{R}(D_1+\eta_{D_1}^{-O(1)})$ and rank at least $R(D_2)$, and repeat the whole process starting with $\mathcal{B}_2$.

Since there is an energy increase of at least $\delta^2$ at every step, this procedure must terminate in at most $M \leq \delta^{-2}$ steps. The result is a $(\delta, \omega, R)$-quadratic regularity partition $\mathcal{B}_M$ whose complexity $D_M$ satisfies the recurrence $D_M \leq \phi_{R}(D_{M-1} + \eta_{D_{M-1}}^{-O(1)})$. Thus, we have an upper bound $C_{qarl}(\delta, \omega, R) \leq D_M$.

What remains is to show that this recurrence relation leads to at most tower-type growth. The following lemma establishes a bound on $\phi_{R}$ by following the proof of \cite[Lemma 3.11]{montreal} (stated as Lemma \ref{lemma:high-rank-refinement-qdr} above). Note that $R(D) = O(D)$ with $R(D) \geq D$ is expected in applications as a consequence of Lemma \ref{lemma:atom-size}.

\begin{lemma}
    \label{lemma:quadratic-rank-refinement-bound}
    Let $R: \mathbb{N} \rightarrow \mathbb{N}$ be a non-decreasing function such that $R(D) \geq D$. Then $\phi_{R}(D) \leq (DR)^\bracketed{D-1}(D)$, where $(DR)^\bracketed{D-1}$ denotes the function $DR$ iterated $D-1$ times, i.e.~$(DR)(x) = D \cdot R(x)$ and $(DR)^{\bracketed{i+1}}(x) = DR((DR)^{\bracketed{i}}(x))$. In particular, if $R(D)$ is of the order $O(D)$, then $\phi_{R}(D) \leq 2^{2^{O(D)}}$.
\end{lemma}
\begin{proof}
    With notation as in Lemma \ref{lemma:high-rank-refinement-qdr}, let the linear layer $\mathcal{B}[1]$ of $\mathcal{B}$ be defined by linear polynomials $L_1, \ldots, L_{D_1}$. Additionally, let the \emph{quadratic layer} $\mathcal{B}[2]$ be defined by quadratic polynomials $Q_1, \ldots, Q_{D_2}$. Here, the quadratic layer of $\mathcal{B}$ is defined analogously to the linear layer, i.e.~as the quadratic factor obtained by restricting the definition of $\mathcal{B}$ to polynomials of degree strictly 2. As such, $\mathcal{B}$ is defined by $(L_1, \ldots, L_{D_1}, Q_1, \ldots, Q_{D_2})$ so that the complexity $D$ of $\mathcal{B}$ satisfies $D = D_1+D_2$.

    Additionally, let $M_i$ denote the unique $n \times n$ matrix such that $L_{Q_i}(x) = Q_i(x) - x^T M_i x$ is a linear polynomial. Recall that $\rank(\mathcal{B}) < R(D)$ if and only if there is a non-zero $\lambda \in \vecsp{p}{D_2}$ such that the matrix $M_{\lambda} = \sum_{i=1}^{D_2} \lambda_{i} M_i$ has rank at most $R(D)$ (see the definition of rank preceding Lemma \ref{lemma:atom-size}). This means that if $\rank(\mathcal{B}) < R(D)$, then any basis $\set{b_1, \ldots, b_r}$ of $\Img{M_{\lambda_{i}}}$ must have size at most $R(D)$.
    
    Without loss of generality, assume that $\lambda_1 \neq 0$, and let $L'_1, \ldots, L'_r$ be linear polynomials defined by $L'_i(x) = b_i^T x$. Then for every $x$, there exists a $\mu \in \vecsp{p}{r}$ such that $M_1 x = \sum_{i=1}^{r} \mu_{i} L'_i(x) + \sum_{i=2}^{D_2} \lambda_{i} M_i x,$ so that 
    $$Q_1(x) = x^T M_1 x + L_{Q_1} = x^T\sum_{i=1}^{r} \mu_{i} L'_i(x) + L_{Q_1}(x) - \sum_{i=2}^{D_2} \lambda_{i} L_{Q_i}(x) + \sum_{i=2}^{D_2} \lambda_{i} Q_i(x).$$
    In particular, if we replace $Q_1$ with the linear polynomials $L'_1, \ldots, L'_r, L_{Q_1} - \sum_{i=1}^{D_2} \lambda_{i} L_{Q_i}$ in the definition of $\mathcal{B}$, the resulting factor $\mathcal{B}^{\bracketed{2}}$ is a refinement of $\mathcal{B}$. This new factor has complexity $D^\bracketed{2} \leq D+R(D) \leq 2R(D)$ and, moreover, $\mathcal{B}^{\bracketed{2}}[2]$ is defined by at most $D_2-1$ quadratic polynomials.

    If the rank of $\mathcal{B}^{\bracketed{2}}$ is not sufficiently high, i.e.~the rank is less than $R(D^{\bracketed{2}}) \leq R(2R(D))$, repeat the argument with $\mathcal{B}^{\bracketed{2}}$ in place of $\mathcal{B}$. The result is a quadratic factor $\mathcal{B}^{\bracketed{3}}$ refining $\mathcal{B}^{\bracketed{2}}$ such that its complexity $D^{\bracketed{3}}$ satisfies
    $$D^{\bracketed{3}} \leq D^\bracketed{2} + R(D^\bracketed{2}) \leq 2R(D) + R(2R(D)) \leq 3R(3R(D))$$
    and $\mathcal{B}^{\bracketed{3}}[2]$ is defined by at most $D_2-2$ quadratic polynomials. Likewise, after another step of the iteration, $\mathcal{B}^{\bracketed{4}}$ will have complexity
    $$D^{\bracketed{4}} \leq 3R(3R(D)) + R(3R(3R(D))) \leq 4R(3R(3R(D))) \leq 4R(4R(4R(D))),$$
    with $\mathcal{B}^{\bracketed{4}}[2]$ defined by at most $D_2-3$ quadratic polynomials and so on.

    Evidently, this process cannot continue for more than $D_2 \leq D$ steps before there are no quadratic polynomials left. As a result, there is some $M \leq D$ such that $\mathcal{B}' = \mathcal{B}^{\bracketed{M}}$ has rank at least $R(D^{\bracketed{M}})$ and complexity $D' = D^{\bracketed{M}} \leq (MR)^{\bracketed{M-1}}(D) \leq (DR)^{\bracketed{D-1}}(D)$, as required. Finally, if $R(D) = O(D)$, then $DR(t) \leq O(t^{2})$ whenever $t \geq D$. Therefore, $\phi_{R}(D) \leq (DR)^{\bracketed{D-1}}(D) \leq D^{O(2^{D-1})} = 2^{2^{O(D)}}$.
\end{proof}

\noindent As $\omega(t)$ is typically of the order $2^{O(t)}$ in applications, consider $\eta_t^{-1} = 2^{O(t)}$ so that
$$D_{i} \leq \phi_{R}\big(D_{i-1} + \omega(D_i)^{O(1)}\big) \leq \phi_{R}\big(D_{i-1}+2^{O(D_i)}\big) \leq \phi_{R}\big(2^{O(D_{i-1})}\big) \leq 2^{2^{2^{O(D_i)}}}.$$
Starting with $D_0 = 0$ and continuing for $M \leq \delta^{-2}$ steps gives an upper bound on $D_M$ that grows like a tower of height $3\ceil{\delta^{-2}}$, as required. This is comparable with the growth of complexity in the linear case $s=1$ (Theorem \ref{theorem:arl}) which is at most a tower of height $\ceil{\epsilon^{-3}}$ \cite{tower-type}.

\if{0}
{\color{teal}
Note that similar bounds may be derived for higher-order arithmetic regularity lemmas. There are two differences compared to the quadratic case. Firstly, the best-known bounds in the inverse theorem (Theorem \ref{theorem:gowers-inverse}) for the $U^{s+1}$-norm are currently exponential for $s \geq 3$ \cite{effective-bounds} \cmnt{although see \url{https://arxiv.org/pdf/2402.17994} and \url{https://arxiv.org/abs/2410.08966}}. All this means, however, is that the recurrence one obtains on $D_i$ has the form $D_i \leq \phi_{R}\big(\exp^{O(1)}(O(D_{i-1}))\big)$, where $\exp^{O(1)}(\cdot)$ denotes $\exp(\cdot)$ iterated $O(1)$ times, i.e.~a tower of fixed height, much like in the quadratic setting.

Secondly, the notion of rank used for arithmetic regularity lemmas of order greater than $2$ has a somewhat different form (see \cite[Definition 1.5]{rank-def}). Nevertheless, a higher-order analogue of Lemma \ref{lemma:high-rank-refinement-qdr}, such as \cite[Lemma 3.18]{decomposition-lemmas}, has a very similar proof. \cmnt{but we need to iterate more times now - will we still be tower?} Considerations on the desired size of atoms (which may be handled by a higher-order version of Lemma \ref{lemma:atom-size}, such as \cite[Lemma 3.2]{trans-invariant}) along with the current best-known bounds on rank \cite{moshkovitz-zhu} still give $R(x) = O(x^2)$. As a result, we get a bound that is a much higher tower, perhaps a tower iterated a fixed number of times, but a tower nonetheless.}
\fi

\end{document}